\documentclass[reqno,11pt]{article}
\usepackage{amsmath, amsthm, amscd,amsfonts, amssymb,hyperref,floatrow}
\usepackage{graphicx, color,dsfont,xcolor}

\numberwithin{equation}{section}

\newtheorem{theorem}{Theorem}[section]
\newtheorem{lemma}[theorem]{Lemma}
\newtheorem{proposition}[theorem]{Proposition}

\newtheorem{assumption}[theorem]{Assumption}

\newtheorem{conjecture}[theorem]{Conjecture}

\renewcommand{\tilde}{\widetilde}          
\DeclareMathSymbol{\leqslant}{\mathalpha}{AMSa}{"36} 
\DeclareMathSymbol{\geqslant}{\mathalpha}{AMSa}{"3E} 
\DeclareMathSymbol{\eset}{\mathalpha}{AMSb}{"3F}     
\renewcommand{\leq}{\;\leqslant\;}                   
\renewcommand{\geq}{\;\geqslant\;}                   

\def \C{ \mathbb  C }
\newcommand{\R}{\mathbb{R}}

\newcommand{\N}{\mathbb{N}}

\newcommand{\ind}{\mathds{1}}
\def \P{ \mathbb P  }
\def \E{ \mathbb E  }
\newcommand \be  {\begin{equation}}
\newcommand \bea {\begin{eqnarray} \nonumber }
\newcommand \ee  {\end{equation}}

\newcommand \ba  {\begin{align}}
\newcommand \ea  {\end{align}}

 \hypersetup{
    linktoc=page,
    linkcolor=blue,          
    citecolor=red,        
    filecolor=blue,      
    urlcolor=cyan,
    colorlinks=true           
}

\definecolor{remi}{rgb}{0,0,0}
\usepackage{titlesec}
    \titleformat{\section}[hang]
        {\color{remi}{}\bfseries\filcenter\large}
        {\thesection.}
        {0.4em}
        {}[]

\title{Marchenko Pastur type theorem for independent MRW processes: convergence of the empirical spectral measure}
\author{}

\begin{document}
\maketitle
\begin{center}
{Romain Allez, R\'emi Rhodes, Vincent Vargas \\
\footnotesize 
 
 CNRS, UMR 7534, F-75016 Paris, France \\
 Universit{\'e} Paris-Dauphine, Ceremade, F-75016 Paris, France} \\

{\footnotesize \noindent e-mail: \texttt{allez@ceremade.dauphine.fr, rhodes@ceremade.dauphine.fr, vargas@ceremade.dauphine.fr}}
\end{center}

\begin{abstract}
We study the asymptotic of the spectral distribution for large empirical covariance matrices composed of independent Multifractal Random Walk processes. The asymptotic is taken as the observation lag shrinks to $0$. In this setting, we show that there exists a limiting spectral distribution whose Stieltjes transform is uniquely characterized by equations which we specify. We also illustrate our results by numerical simulations.  
\end{abstract}

\footnotesize



\normalsize

\tableofcontents

\section{Introduction}
Since the seminal work of Mar\~cenko and Pastur \cite{MP} in 1967, there has been growing interest in studying the asymptotic of large empirical covariance matrices. These studies have found applications in many fields 
of science: physics, telecommunications, information theory and finance, etc... The main motivation of this work stems from finance: the study of covariance matrices is a crucial tool for minimizing 
the risk $\mathcal{R}_w$ of a portfolio $w$ that invests $w_i$ in asset number $i$. Indeed, if we denote by $r_i$ the price variation of asset $i$, $\mathcal{R}_w$ can be defined as 
the variance of the random variable 
$\sum_i w_i r_i$ and can be computed in terms of the covariance matrix $\bf R$ of the $r_i$ (defined as ${\bf R}_{ij}=\E[r_ir_j]$): 
\begin{equation*}
\mathcal{R}_w =  w^t {\bf R} w\,.
\end{equation*}
Of course, practitioners do not have access to $\bf R$; instead, they must consider a noisy empirical estimator of $\bf R$, which consists of a large empirical covariance matrix. A key tool in distinguishing noise from real correlations is the study of the eigenvalues of the empirical covariance matrix: we refer to  \cite{JPBshortreview}, \cite{JPBOldLaces} for more extended discussions on the applications of large empirical covariance matrices in finance and in particular in portfolio theory.   


We will work in a high frequency setting: we consider $N$ stock price processes $X_i(t)$ for $i=1,\dots,N$ that evolve continuously 
with respect to time $t\in [0;1]$   
but we observe those prices {\it only on a discrete finite grid} $\{j/T,j=1,\dots,T\}$ where $T$ is the number of observations. 
Using this discrete grid, we can compute the price variations $r_i(j)$ (that we will abusively call {\it returns}) for each asset price $X_i$ on every time interval $[(j-1)/T;j/T]$ by:   
\begin{equation*}
r_i(j) := X_i(\frac{j}{T}) - X_i(\frac{j-1}{T}).
\end{equation*}
Then, we define the $N\times T$ matrix $X_N$ such that $X_N(ij)= r_i(j)$ that enables to define the empirical covariance matrix $R_N$ as follows $$R_N := X_N X_N^t\,.$$
In this work, we will be interested in the statistics of the symmetric matrix $R_N$ and in particular in its spectrum, or more precisely, in its  limiting spectral distribution in the limit of large matrices 
(i.e. when $N\rightarrow \infty$) for different models
of the {\it i.i.d.} random continuous processes $(X_i(t)),i\in \{1,\dots,N\}$ (see below for precise definitions).
For this purpose, the Mar\~cenko-Pastur paper enables to deal with the case where stock prices follow independent Brownian motions. 
More precisely, in this case, the matrix $X_N$ is defined as:
\be
X_N(ij) = B_i\left(\frac{j}{T}\right) -B_i\left(\frac{j-1}{T}\right)  
\ee
where the $B_i$ are {\it i.i.d.} standard Brownian motions. 

If $\lambda_1, \dots, \lambda_N$ are the eigenvalues of $R_N$, 
the empirical spectral distribution of the matrix $R_N$ is the probability measure defined by: 
\be
\mu_{R_N} = \frac1N \sum_{i=1}^{N} \delta_{\lambda_i}.
\ee
The Mar\~cenko-Pastur (MP) result states that, in the limit of large matrices $N,T\to\infty$ with $N/T \to q \in (0,1]$, the empirical spectral distribution $\mu_{R_N}$ 
weakly converges (almost surely) to a probability measure whose density $\rho(x)$ is: 
\be
\rho(x) = \frac{1}{2\pi q} \frac{\sqrt{(\gamma_+ - x)(x - \gamma_-)}}{x} \mathds{1}_{[\gamma_-,\gamma_+]} dx
\ee
where $\gamma_\pm = 1+q \pm 2\sqrt{q}$.

Independently of the aforementioned work on random matrix theory, much work has been devoted to studying the statistics of financial stocks. It turns out that most financial assets (stocks, indices, etc...) 
possess universal features, called stylized facts. In short, one can observe empirically the following properties (the list below is obviously non exhaustive) for asset returns on financial markets: 
\begin{itemize}
\item The returns are multifractal; in particular on short scales, they are heavy tailed but tend to have distribution closer to the Gaussian law on larger scales. 
\item The volatility fluctuates randomly and follows approximately a lognormal distribution.
\item While the returns are rapidly decorrelated, the volatility exhibits long range correlations following a power law. 
\end{itemize}
We refer to the references \cite{BoPobook,Cont} for a discussion on this topic. Many models have been proposed in the literature that take into account 
these stylized facts. Among them, there has been growing interest in the lognormal Multifractal Random Walk (MRW) model introduced in \cite{Bacry} (see also \cite{Bacry2,Vincent_vol}). The lognormal MRW model satisfies several of the so-called stylized facts, 
but a few of them remain unchecked such as asymmetry of returns and Leverage effect (see \cite{BoMaPo}). The lognormal MRW is simply defined as: 
\be
X(t) = B\left(M[0,t]\right)
\ee 
where $B$ is a standard Brownian motion and $M$ is an independent lognormal multifractal random measure (MRM for short) formally defined, for $t\geq 0$, by: 
\begin{equation*}
M[0;t] = \int_0^t e^{\omega(x)-\frac{1}{2} \E[\omega(x)^2]} {\rm d}x\,,
\end{equation*}
where $(\omega(x))_{x\in \R}$ is a "gaussian field" whose covariance kernel $K$ is 
\begin{equation*}
K(x,y) = \gamma^2 \ln_+\left(\frac{\tau}{|t-s|}\right)\,,
\end{equation*} 
where $\ln_+ x=\text{max}(\ln x, 0)$. The two parameters $\gamma^2$ and $\tau$ are respectively called intermittency parameter and integral scale (or correlation length) of the lognormal random multifractal measure $M$. 

\begin{figure}[h!btp] 
	\center
         \includegraphics[scale=0.8]{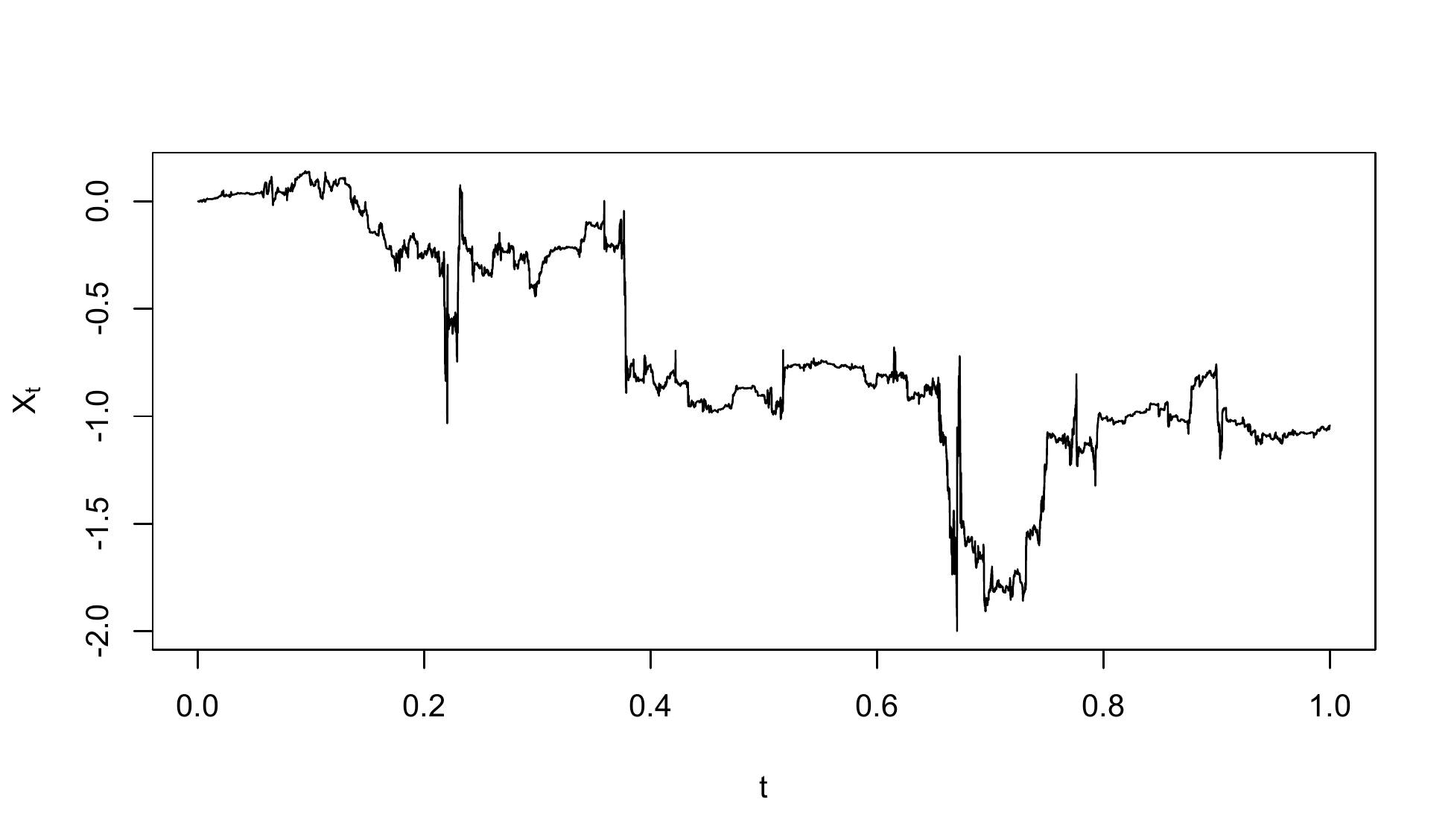}
         \caption{Simulated path of a multifractal random walk with intermittency parameter $\gamma^2=1$ and with integral scale $\tau=1/4$. Note the intermittent bursts in volatility.}\label{fig1_mrw}
\end{figure}

Fig. \ref{fig1_mrw} represents a simulated path of a lognormal MRW $X(t)=B(M([0;t]))$ where $B$ is a standard Brownian motion independent of the multifractal random measure $M$ with intermittency parameter 
$\gamma^2=1$ and integral scale $\tau=1/4$. The reader can find a more precise reminder of the construction/definition of a more general class of Multifractal Random Measure (MRM), as well as (standard) notations 
used throughout the paper in section \ref{cons:loglevy}.

We thus aim at studying the large sample covariance matrices where the underlying price processes evolve as lognormal MRW. 
More precisely, the matrix $X_N$ is defined, for $1\leq i \leq N, 1 \leq j \leq T$, as: 
\begin{equation}\label{entries}
X_N(ij) =   B^i({M^i(0,\frac{j}{T})})- B^i({M^i(0,\frac{j-1}{T})})
\end{equation}
where the $B_i$ are {\it i.i.d.} Brownian motions and the $M_i$ are {\it i.i.d.} lognormal MRM independent of the $B_i$. Let us mention the work \cite{LiZheng} which considers high frequency covariance matrices in the context of diffusion processes (see also \cite{Ros} for studies of high frequency large empirical covariance matrices motivated by financial applications). The processes described by (\ref{entries}) are typically not diffusions.

In the spirit of the MP Theorem, the purpose of this work is to characterize the limit of the empirical spectral measure 
$\mu_{R_N}$ when $N, T \to \infty$ with $N/T \to q \in (0,1]$. It is interesting to understand how the {\it long-memory} volatility   
process affects the covariance matrix in the limit of large matrices. In particular, we will see that the intermittent volatility has the effect to spread the spectrum of the covariance matrix $R_N$ in a wider region of $\R_+$.  
Indeed the spectral density has a compact support $[\gamma_-;\gamma_+]$ in the Mar\~cenko-Pastur setting (in which the prices follow Brownian motions) whereas it has an infinite support with a tail that gets heavier 
as the intermittency parameter grows. We mention that our results can be extended to many different auto-correlated volatility processes. 

The effect of the integral scale $\tau$ on the empirical covariance matrix $R_N$ is also very interesting in the context of price variations measured on a very short scale (high frequency). The high frequency 
case corresponds to large 
values of the parameter $\tau$ while low frequency case corresponds to small values of $\tau$. Indeed, if $X$ is a lognormal MRW with integral scale $\tau$, then the process 
$\widehat{X}(t)$ defined on $[0;1]$ as $\widehat{X}(t)= X(t/2)$ is a lognormal MRW with integral scale $2\tau$. 
Note that this discussion on high freqency measurement is irrelevant in the MP case when asset prices follow independent Brownian motions
since, in this model, the distribution of price variations is the same on any scale: it is Gaussian, only the variance will change with the scale and up to the variance parameter the limiting spectral distribution
will always be the same at different scales. However, if asset prices follow lognormal MRW (or even another process with a correlated in time volatility process), the price variations measured on small scales will have a 
distribution with higher kurtosis (i.e. the probability mass of the tail is heavier)
and therefore the spectrum of the empirical covariance matrix $R_N$ should be affected by decreasing the measurement scale. We therefore expect stronger right tail for the spectral distribution. 
The numerical analysis of our results indeed confirms this guess: the larger the integral scale is, the heavier is the right tail.   

Here, we are mainly interested in the case where asset prices follow lognormal MRW but we will also present our results for two other related models 
where asset prices follow independent Brownian motions with a time change, which can be thought of as a volatility 
process with memory (i.e. the volatility process is correlated in time). 

%
%
%


The next sections are organized as follows. In section 2, we remind the definition of MRW and introduce the main notations of the paper. In section 3, we state our main theorems which are characterizations of the 
limiting spectral measure of $R_N$ through its Stieltjes transform for different types of underlying processes $X$. These equations are tedious to invert analytically and it is hard to extract the properties (continuity,
tails of the distribution) of the associated spectral 
density. 
In section 4, we invert these equations numerically so as to get informations on the spectral measure of the covariance matrix $R_N$ as $N\rightarrow \infty$ and we check the validity and applicability of our results 
using numerical simulations. 
The proofs appear in section 5 with some auxiliary lemmas proved in the appendix. The strategy of our proofs is classical among the random matrix literature 
(the so-called resolvent method) as it relies on the Schur recursion formula for the Stieltjes transform; in particular, we follow the approach of \cite{BenArousGuionnet}. The main difficulty lies 
in handling the Stieltjes transforms in a multifractal setting.

\section{Background, notations and main results}

\subsection{Reminder of the construction of  MRM}\label{cons:loglevy}
To fix precisely the notations that we will use throughout the paper, we quickly remind the main steps of the construction of Multifractal Random Measures (MRM). The description is necessarily concise and the reader is referred to \cite{Bacry} for further details. In particular, we use the same notations as in \cite{Bacry} to facilitate the reading.   We consider the characteristic function of an infinitely divisible random variable $Z$, which can be written as $\E[e^{ipZ}]=e^{\varphi(p)}$ where (L\'evy-Khintchine's formula):
\be\label{char}
\varphi(p)=i mp-\frac{1}{2}\gamma^2p^2+\int_{\R^*} (e^{ipx} -1)\,\nu(dx)
\ee
and $\nu(dx)$ is a so-called L\'evy measure (ie satisfying $\int_{\R^*}\min(1,x^2)\,\nu(dx)<+\infty$) together with the following additional assumption:
\be
\int_{[-1,1]}|x|\,\nu(dx)<+\infty,
\ee
so that its characteristic function perfectly makes sense as written in \eqref{char}.
We also introduce the Laplace exponent $\psi$ of $Z$ by $\psi(p)=\varphi(-ip)$ for each $p$ such that both terms of the equality make sense, and we assume that the following renormalization condition holds: $\psi(1)=0$.  

We further consider the half-space $ S=\{(t,y);t\in\R,y\in\R^*_+\}$, with which we associate the measure (on the Borel $\sigma$-algebra $\mathcal{B}(S)$):
\be\label{theta}
\theta(dt,dy)=y^{-2}dt\,dy. 
\ee
Then we consider an independently scattered infinitely divisible random measure $\mu$ associated to $(\varphi,\theta)$ and distributed on $S$. 
 
Then we define a process $\omega_\epsilon$ for $\epsilon>0$ by the following.
Given a positive parameter $\tau$, let us  define the function $f:\R_+\to \R$ by:
$$ f(r)=\left\{
\begin{array}{ll}
 r, &  \text{ if } r\leq \tau  \\
 \tau &  \text{ if } r\geq \tau
\end{array}\right.\,.$$ 
The cone-like subset $A_\epsilon(t)$ of $S$ is defined by:
\begin{equation}
A_\epsilon(t)=\{(s,y)\in S; y\geq \epsilon, -f(y)/2\leq s-t\leq f(y)/2\}.
\end{equation}
We then define the stationary process $(\omega_\epsilon(t))_{t \in \R}$ by:
\begin{equation}
\omega_\epsilon(t) = \mu\left(A_\epsilon(t)\right).
\end{equation}

The Radon measure $M$ is then defined as the almost sure limit (in the sense of weak convergence of Radon measures) by: 
$$M(A)=\lim_{\epsilon\to 0^+} M_\epsilon(A)=\lim_{\epsilon\to 0^+} \int_Ae^{\omega_\epsilon(r)}\,dr$$
for any Lebesgue measurable subset $A\subset \R $. The convergence is ensured by the fact that the family $(M_\epsilon(A))_{\epsilon>0}$ is a right-continuous positive martingale. The structure exponent of $M$ is defined by:
$$\forall p \geq 0, \quad \zeta(p)= p-\psi(p) $$ for all $p$ such that the right-hand side makes sense.
The measure $M$ is different from $0$ if and only if there exists $\epsilon>0$ such that $ \zeta(1+\epsilon)>1$, (or equivalently $\psi'(1)<1$). In that case, we have:

\begin{theorem}\label{prop:existMRM}
The measure $M$ is stationary and satisfies the  {\bf exact stochastic scale invariance property}: for any $\lambda\in ]0,1]$, 
$$ (M(\lambda A))_{A\subset B(0,\tau)}\stackrel{{\rm law}}{=}(\lambda  e^{\Omega_\lambda}  M(A))_{A\subset B(0,\tau)},$$ where $\Omega_\lambda$ is an infinitely divisible random variable, independent 
of $ (M(A))_{A\subset B(0,T)}$, the law of which is characterized by:
$$ \E[e^{ip \Omega_\lambda}]=\lambda^{-\varphi(p)}.$$
\end{theorem}

\subsection{Notations}

Let $N$ and $T:=T(N)$ be two integers, the aim of this paper is to compute the empirical spectral measure of the matrix $R_N :=  X_N  {}^{t} \! X_N$ as $N \to \infty$, 
where $X_N$ is a $N\times T$ real matrix the entries of which are given by \eqref{entries}. Recall that the number $N$ of sampled processes is supposed to be comparable with the sample size $T:=T(N)$, and more precisely, 
we will suppose in the following that there exists a parameter $q\in]0,1]$ such that:
\begin{equation}\label{regime}
\lim_{N\to \infty}\frac{N}{T}=q.
\end{equation}
We further set $\widetilde{R}_N :=  {}^{t} \! X_N  X_N$, and if $M$ is a symmetric real matrix, we will denote by $\mu_M$ the empirical spectral measure of $M$.  

Define the $(T+N) \times (T+N)$ matrix $B_N$ by:
\begin{equation*}
B_N =
 \begin{pmatrix}
         0 &  {}^{t} \! X_{N} \\
         X_{N} & 0
\end{pmatrix}.
\end{equation*}

We also define for $z \in \C \setminus \R$,
 \begin{equation*}
A_N(z)= \left(zI_{T+N} - B_N \right) =
 \begin{pmatrix}
         z I_{T} & -  {}^{t} \! X_{N} \\
         -X_{N} & z I_{N }
\end{pmatrix}.
\end{equation*}
Notice that $$B_N^2=\begin{pmatrix}
         \widetilde{R}_N&0 \\
       0&  R_{N}
\end{pmatrix}$$ and that the eigenvalues of $\widetilde{R}_N$ are those of $R_N$ augmented with $T-N$ zero eigenvalues. 
We thus have:
\begin{equation}\label{rel:rn}
\mu_{B_N^2}=2\frac{N}{N+T}\mu_{R_N}+\frac{T-N}{N+T}\delta_0, 
\end{equation}
where $\delta_x$ stands for the Dirac mass at $x$.
Combining this equality with the relation
\begin{equation}
\int f(x)\mu_{B_N^2}(dx)=\int f(x^2)\mu_{B_N}(dx) 
\end{equation}
true for all bounded continuous functions $f$ on $\R$, we see that it is sufficient to study the weak convergence of the spectral measure of $B_N$ for the study of the convergence 
of the spectral measure $\mu_{R_N}$.  

We will thus work on the (weak) convergence of the spectral measures $\mu_{B_N}$ and $\E\left[\mu_{B_N}\right]$ in the following. To that purpose, it is sufficient to prove the convergence of the Stieltjes transform 
of these two measures. 
Recall that, for a probability measure $\mu$ on $\R$, the Stieltjes transform $G_{\mu}$ of $\mu$ is defined, for all $z \in \C \setminus \R$, as:  
\begin{equation}
G_{\mu}(z)= \int_\R\frac{ 1}{z-x}\mu(dx).
\end{equation}
and one can note that:
\begin{equation}\label{rel:trace}
G_{\mu_{B_N}}(z)= \frac{1}{N+T}{\rm Trace}(G_N(z)),
\end{equation}
where we have set:
\begin{equation}
G_N(z) = \left( A_N(z) \right)^{-1}.
\end{equation}
Hence, we have to investigate the convergence of the right-hand side of \eqref{rel:trace}. Let us introduce the two following complex measures $L^{1,z}_N$ and $L^{2,z}_N$ such that, for all bounded and 
measurable function $f:[0,1] \to \R$:
\begin{align*}
L^{1,z}_N(f) &= \frac{1}{T} \sum_{k=1}^{T} f\left(\frac{k}{T}\right) G_N(z)_{kk}  \\
L^{2,z}_N(f) &= \frac{1}{N} \sum_{k=1}^{N}  f\left(\frac{k}{N}\right) G_N(z)_{k+T,k+T}
\end{align*}
Clearly, we have the relation
\begin{equation}\label{rel:trace2}
\frac{1}{N+T}{\rm Trace}(G_N(z))=\frac{T}{N+T}L^{1,z}_N([0,1])+\frac{N}{N+T} L^{2,z}_N([0,1])
\end{equation}
so that it suffices to establish the convergence of the two complex measures $L^{1,z}_N$ and $L^{2,z}_N$. 

\section{Main results}
\subsection{Lognormal multifractal random walk}\label{lognormal_multifractal}
We first present our results when the process $X(t)$ is a lognormal multifractal random walk, i.e. $X(t) = B(M[0;t])$ 
where $M$ is the MRM whose characteristic and structure exponent (see section \ref{cons:loglevy}) are respectively given by:
\begin{align*}
\varphi(q) &= -i\frac{\gamma^2}{2}q - \frac{\gamma^2}{2} q^2, \\
\zeta(q) &= (1+\frac{\gamma^2}{2}) q - \frac{\gamma^2}{2} q^2.
\end{align*}

We will make the assumption that the intermittency parameter $\gamma^2$ is small enough so as to overcome in our proofs the strong correlations of the model. 
\begin{assumption} \label{cond:gamma^2} More precisely, let us suppose that:
\begin{equation}
\gamma^2 < \frac{1}{3}.
\end{equation}
\end{assumption}
Though we conjecture that our results hold as soon as the measure $M$ is non degenerated, i.e. $\gamma^2 < 2$ (see \cite{Bacry}), Assumption \ref{cond:gamma^2} is largely
sufficient to cover most practical applications. For instance, in financial applications or in the field of turbulence, $\gamma^2$ is found empirically around $2.10^{-2}$.

We can now state our main result about the convergence of the empirical spectral measures and mean empirical spectral measures of the matrices $B_N$ and $R_N$: 
\begin{theorem}\label{main2}
i) There exists a probability measure $\upsilon$ on $\R$ such that the two mean spectral measures $\E[\mu_{B_N}]$ and $\E[\mu_{R_N}]$ converge weakly respectively towards
the two probability measures
$\frac{2q}{1+q}\upsilon+\frac{1-q}{1+q}\delta_0$ and $\upsilon\circ (x^2)^{-1}$ as $N$ goes to $\infty$, where $\upsilon\circ (x^2)^{-1}$ is the push-forward of the measure $\upsilon$ by the mapping $x\mapsto x^2$.

ii) The two spectral measures $\mu_{B_N}$ and $\mu_{R_N}$ converge weakly in probability respectively to the two probability measures $\frac{2q}{1+q}\upsilon+\frac{1-q}{1+q}\delta_0$ and $\upsilon\circ (x^2)^{-1}$
as $N$ goes to $\infty$. More precisely, for any bounded and continuous function $f$, $\int f(x) \mu_{R_N}(dx)$ converges in probability to $\int f(x) \upsilon\circ (x^2)^{-1}(dx)$. 

iii) Let $N_k$ be an increasing sequence of integers such that $\sum_{k=1}^{\infty} N_k^{-1} < +\infty$, then the two sequences $\mu_{B_{N_k}}$ and $\mu_{R_{N_k}}$ converge weakly almost surely to 
the two probability measures
$\frac{2q}{1+q}\upsilon+\frac{1-q}{1+q}\delta_0$ and $\upsilon\circ (x^2)^{-1}$ as $k$ goes to $\infty$. 
\end{theorem}

Theorem \ref{main2} is implied by \eqref{rel:trace}, \eqref{rel:trace2} and by Theorem \ref{main1}:
\begin{theorem}\label{main1}
i) The measures $\E[L^{1,z}_N]$ and $\E[L^{2,z}_N]$ converge weakly towards two complex measures. More precisely, there exist a unique $\mu^2_z\in\C$ and a unique bounded measurable 
function $K_z(x)$ over $[0,1]$ such that, for all bounded and continuous function $f$ on $[0,1]$, we have respectively: 
\begin{align*}
\E\left[L^{1,z}_N(f)\right]\to_{N\to \infty} \int_0^1K_z(x)f(x)\,dx,\\
\E\left[L^{2,z}_N(f)\right]\to_{N\to \infty} \mu^2_z\int_0^1f(x)\,dx.
\end{align*}

ii) In addition, we have the following relation between $\mu^2_z\in\C$ and $K_z(x)$: 
\begin{equation}
\int_0^1K_z(x)\,dx=q\mu^2_z+\frac{1-q}{z}
\end{equation}

iii) Furthermore, there exists a unique probability measure $\upsilon$ on $\R$ whose Stieltjes transform is $\mu^2_z$, meaning that for all $z \in\C\setminus\R$,
\begin{equation}\label{funrv}
\mu^2_z=\int_\R\frac{\upsilon(dx)}{z-x}. 
\end{equation}
\end{theorem}

It is important to state a characterization of the probability measure $\upsilon$: it is done by means of its Stieltjes transform $\mu^2_z$:
\begin{theorem}\label{main3}
The constant $\mu_z^2$ and the bounded function $K_z(x)$ are uniquely determined for all $z \in\C\setminus\R$, by the following system of equations:
\begin{align}\label{eqlim:MRM_lognormal}
\mu^2_z &= \E\left[\left(z - \int_0^1  K_z(t) M(dt) \right)^{-1}\right], \\ 
 K_z(x) &= \left( z - q \E\left[ \left(  z - \int_0^1 \Big(\frac{\tau}{|t-x|}\Big)_+^{\gamma^2} K_z(t) M(dt)  \right)^{-1}\right]\right)^{-1} \label{eqlim:lognormalMRM}
\end{align}
where\footnote{The notation $(\cdot)_+$ is a shortcut for $\max(\cdot,1)$.} $M$ is the MRM with structure exponent $\zeta(q)=(1+\gamma^2/2)q - q^2\gamma^2/2$.
\end{theorem}

Let us notice that one can give a precise meaning to (\ref{eqlim:lognormalMRM}) for all $\gamma^2 \in [0,2[$. Indeed, we can define for all $x \in [0,1]$ and all continuous function $f$, the following almost 
sure limit as a definition:
\begin{equation}\label{eq:definitiontendue}
\int_0^1 \Big(\frac{\tau}{|t-x|}\Big)_+^{\gamma^2} f(t)  M(dt) = \lim_{\eta \to 0} \int_{t \in [0,1]; |t-x|> \eta} \Big(\frac{\tau}{|t-x|}\Big)_+^{\gamma^2} f(t)  M(dt)
\end{equation}
Note that the above limit exists almost surely since, for $x$ fixed: 
\begin{equation*}
\frac{\ln M[x-\epsilon_k,x+\epsilon_k]}{\ln \epsilon_k} \underset{k \to \infty}{\rightarrow} 1+\frac{\gamma^2}{2}, \;  a.s.
\end{equation*}
where $\epsilon_k=\frac{1}{2^k}$. One can also check with this definition that we have:   
\begin{equation*}
\int_0^1 \Big(\frac{\tau}{|t-x|}\Big)_+^{\gamma^2} f(t)  M(dt) = \lim_{\epsilon \to 0} \int_{0}^1 e^{cov(\omega_{\epsilon}(t),\omega_{\epsilon}(x))}  f(t)  e^{\omega_{\epsilon}(t)}dt
\end{equation*}

\begin{conjecture}\label{conj}
With this extended definition, we conjecture that theorem \ref{main3} holds in the lognormal multifractal case for all $\gamma^2 \in [0,2[$ and thus that the limiting equations can be obtained by the ones 
of theorem \ref{main3lognorm} (see below) with $2W=\omega_{\epsilon}$ as $\epsilon \to 0$.
\end{conjecture}

\subsection{General multifractal random walk}
We now look at the more general case when the change of time is a measure $M$ for which the function $\varphi(q)$ is given by (\ref{char}) 
and the structure exponent by $\zeta(q) = q -\psi(q)$ with $\psi(q)=\varphi(-iq)$.

We still have to make an assumption to avoid the issue of strong correlations. In this more general setting, Assumption (\ref{cond:gamma^2}) becomes: 

\begin{assumption}\label{ass:MRM} Assume that the structure exponent of the MRM satisfies the condition: 
\begin{equation}\label{condzeta}
\zeta(2)>5-4\zeta'(1).
\end{equation}
and that there exists $\delta > 0$ such that:
\begin{equation}
\zeta(2+\delta) > 1. 
\end{equation}
\end{assumption}
As in the previous section, we conjecture that our results hold as soon as the measure $M$ is non degenerated, i.e. (see \cite{Bacry}) $\zeta(1+\epsilon)>1$ for some $\epsilon>0$. 

Theorems \ref{main2} and \ref{main1} remain unchanged for this more general context. Theorem \ref{main3} becomes:
\begin{theorem}\label{main4}
The constant $\mu_z^2$ and the bounded function $K_z(x)$ are uniquely determined for all $z \in\C\setminus\R$, by the following system of equations:
\begin{align}\label{eqlim:MRM}
\mu^2_z &= \E\left[\left(z - \int_0^1  K_z(t) M(dt) \right)^{-1} \right], \\ 
 K_z(x) &= \left( z - q \E\left[ \left(  z - \int_0^1 \Big(\frac{\tau}{|t-x|}\Big)_+^{\kappa} K_z(t) Q(dt)  \right)^{-1}\right] \right)^{-1}
\end{align}
with $\kappa=\psi(2)$
and where $M$ is the MRM whose characteristic and structure exponent are respectively $\varphi(q),\zeta(q)$ and where the random Radon measure $Q$ is defined, conditionally on $M$, as the almost sure weak 
limit as $\epsilon$ goes to $0$ of the family of random measures $Q_\epsilon(dt):=e^{\overline{\omega}_\epsilon(t)} M(dt)$ where, for each $\epsilon>0$, the random process $\overline{\omega}_\epsilon$ is independent of
$M$ and defined as $\overline{\omega}_\epsilon(t)=\overline{\mu}(A_\epsilon(t))$ where $\overline{\mu}$ is the independently scattered log infinitely divisible random measure 
associated to $(\bar{\varphi},\theta(\cdot \cap A_0(x)))$ with: 
\begin{equation}
\bar{\varphi}(p) = ip (\gamma^2-\kappa) +\int_\R(e^{ipx}-1)(e^x-1)\nu(dx).
\end{equation}
 \end{theorem}

\subsection{Lognormal random walk}

Let us mention that one can easily adapt the methods used to prove the above theorems in the simpler case (lognormal case) where $X(t)$ is defined, for all $t\in [0;1]$, by: 
\be
X(t) = B\left(\int_0^te^{2W(s)} ds \right),
\ee
where $(W(s))_{s\in[0;1]}$ is a stationary gaussian process with expectation $m$ and stationary covariance kernel $k$. The normalization 
will be chosen such that: $m=-k(0)$. 
 
In this context, the entries of $X_N$ are given, for $1\leq i \leq N, 1 \leq j \leq T$ by:
\begin{equation}\label{entriesG}
X_N(ij) =  \frac{1}{\sqrt{T}} e^{W_i(\frac{j}{T})}B^i_j := r_i(j)
\end{equation}
where the $(B^i_j)_{ij}$ are i.i.d standard centered Gaussian random variables and the $W_i$ are i.i.d stationary Gaussian processes with expectation $m$ and stationary covariance kernel $k$. Indeed, if one 
makes the following extra assumption:
\begin{assumption}\label{ass:lognorm} Assume that for some constants $C>0$ and $\beta>0$, the covariance kernel $k$ 
satisfies: 
$$\forall x\in\R,\quad |k(x)-k(0)|\leq C|x|^\beta.$$
\end{assumption}

With the same notations as in the previous section, we can now state the following theorem under assumption \ref{ass:lognorm}:

\begin{theorem}\label{main3lognorm}
The system of equations for $\mu_z^2$ and $K_z(x)$ becomes:
\begin{align}
\mu^2_z &= \E\left[\left(z - \int_0^1  K_z(t) e^{2W(t)}\,dt \right)^{-1}\right] \\ 
K_z(x) &=\left( z - q \E\left[ \left(z - \int_0^1  K_z(t) e^{4k(t-x)} e^{2W(t)} \,dt \right)^{-1} \right]  \right)^{-1}.
\end{align}
where $(W(t))_{t\in[0;1]}$ is a stationary gaussian process with expectation $m$ and stationary covariance kernel $k$.
\end{theorem}

\section{Numerical results and computer simulations}

In this section, we are interested in the case handled in sub-section \ref{lognormal_multifractal}, in which the price of an asset evolves as a lognormal multifractal random walk. 
We want to extract informations on the spectral density $\upsilon\circ (x^2)^{-1}$ 
of the covariance matrix $R_N$ in the limit of large matrices. This section will also give evidence that 
our equations are easy to use in practice for applications. 

The information on the measure $\upsilon$ is entirely contained in its Stieltjes transform $\mu_z^2$ which is the unique solution of the system of equations \eqref{eqlim:MRM_lognormal} and \eqref{eqlim:lognormalMRM}. 
Let us admit for clarity at this point that the measure $\upsilon$ admits a continuous density, at least on the set $\R\setminus\{0\}$. One should be able to show that this is indeed true using the two equations 
\eqref{eqlim:MRM_lognormal} and 
\eqref{eqlim:lognormalMRM} that characterize the probability measure $\upsilon$. Under this continuity assumption for $\upsilon(x)$, we can re-find the density $\upsilon(x)$ from $\mu_z^2$ by the relation 
\begin{equation}\label{eq:recup_res}
\lim_{\epsilon\rightarrow 0} \frac{1}{\pi} \Im(\mu_{x-i\epsilon}^2) = \upsilon(x)\,.
\end{equation}
Note that we just need to find the unique family of functions $(K_z(x))_{x\in [0;1]}$ for $z\in \C\setminus\R$ near the real line, that verifies the fixed point equation \eqref{eqlim:lognormalMRM}. Indeed, 
knowing $(K_z(x))_{x\in [0;1]}$, we can compute $\mu_z^2$ by using equation \eqref{eqlim:MRM_lognormal}, or even simpler, the additional relation that we stated above
\begin{equation}
\int_0^1K_z(x)\,dx=q\mu^2_z+\frac{1-q}{z}\,.
\end{equation}

Let $\mathcal{C}([0;1],\C)$ be the space of bounded functions from $[0;1]$ to $\C$.
For $z\in \C\setminus \R$ fixed, the idea to find $(K_z(x))_{x\in [0;1]}$ is the fixed point method due to Picard. Let us introduce the operator 
$T:\mathcal{C}([0;1],\C) \rightarrow \mathcal{C}([0;1],\C)$ by setting, for $g \in \mathcal{C}([0;1],\C)$ and for all $x\in[0,1]$:
\begin{equation}
Tg(x) = \frac{1}{z- q\E\left[\left( z - \int_0^1 \left(\frac{\tau}{|t-x|}\right)_+^{\gamma^2} g(t) M(dt) \right)^{-1}\right]} \,.
\end{equation}
It can easily be shown (see sub-section \ref{section:unicity}) that if $z\in\C\setminus\R$ is sufficiently far from the real line, then the operator $T$ is contracting and therefore admits a unique fixed point 
$K_z(\cdot)$ in $\mathcal{C}([0;1],\C)$.
To find the fixed point $K_z$, we will iterate the operator $T$ starting from any fixed initial function $K_z^{(0)}$.  We know that, for $z$ such that the operator $T$ is contracting, the $n$-th iteration of the 
function
$K_z^{(n)}:=T(K_z^{(n-1)})$ converges to the unique fixed point $K_z$.
In fact, numerically, there is no need in applying the iteration on $T$ for $z$ such that $T$ is contracting (i.e. for $z$ far from the real line) and one can apply the Picard method directly 
near the real line\footnote{Recall that, in view of equation \eqref{eq:recup_res}, we are interested in the value of the Stieltjes transform near the real line.} and find the fixed point after a reasonable number of 
iterations of the operator $T$.  

The multifractal lognormal random measure $M(dt)$ and multifractal random walk are simulated through the standard method by simulating first, with the use of fast Fourier transform, a gaussian process 
with covariance function given for $\eta >0$ small by  
\begin{equation*}
K_\eta(|t-s|) = \gamma^2 \ln_+(\frac{\tau}{|t-s|+\eta}) \,.
\end{equation*}
The lognormal multifractal random measure and random walk are then constructed from this gaussian process through the standard formulas (see e.g. \cite{Bacry,Vincent_vol}).

\begin{figure}[h!btp] 
	\center
         \includegraphics[scale=0.5]{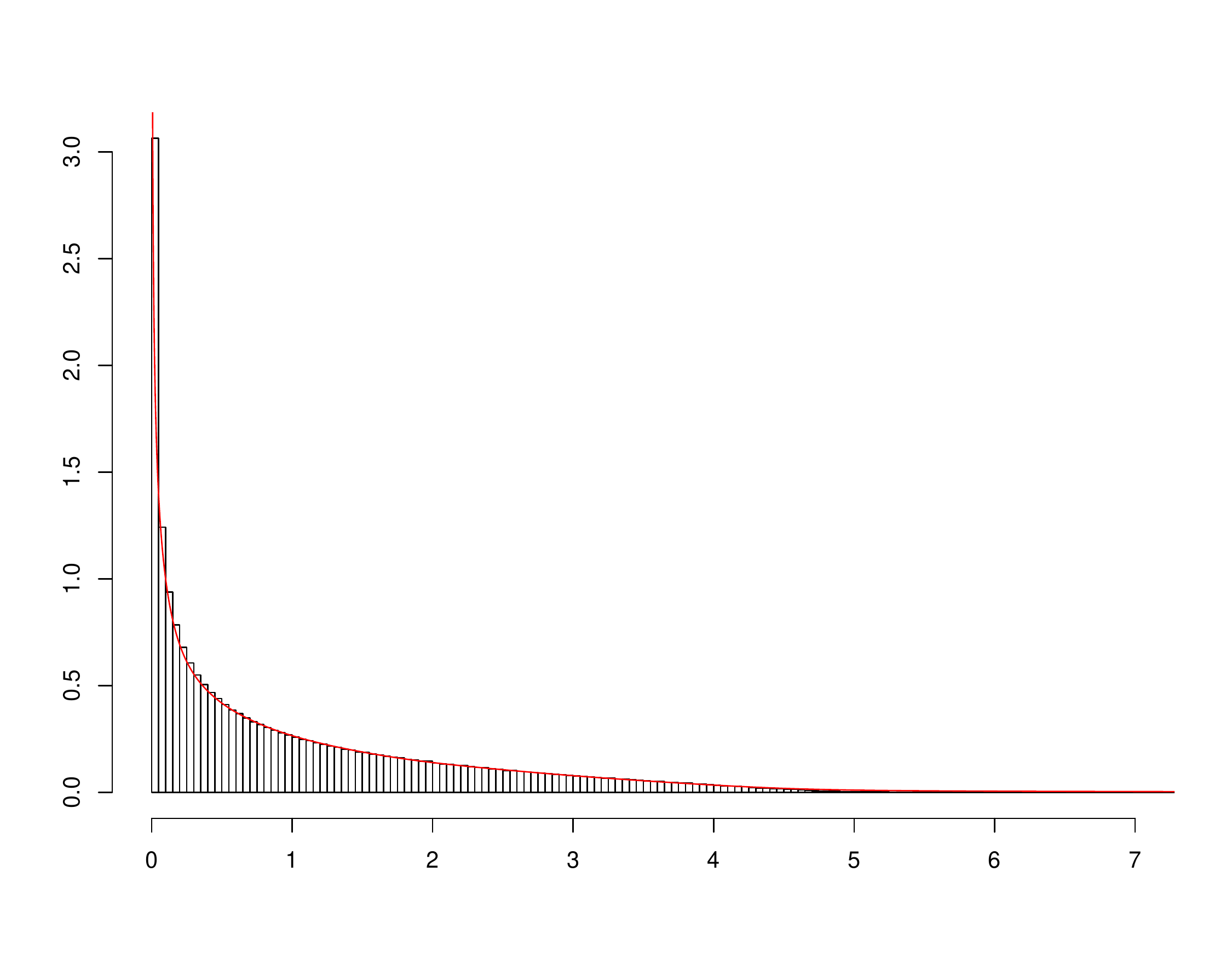}
         \includegraphics[scale=0.5]{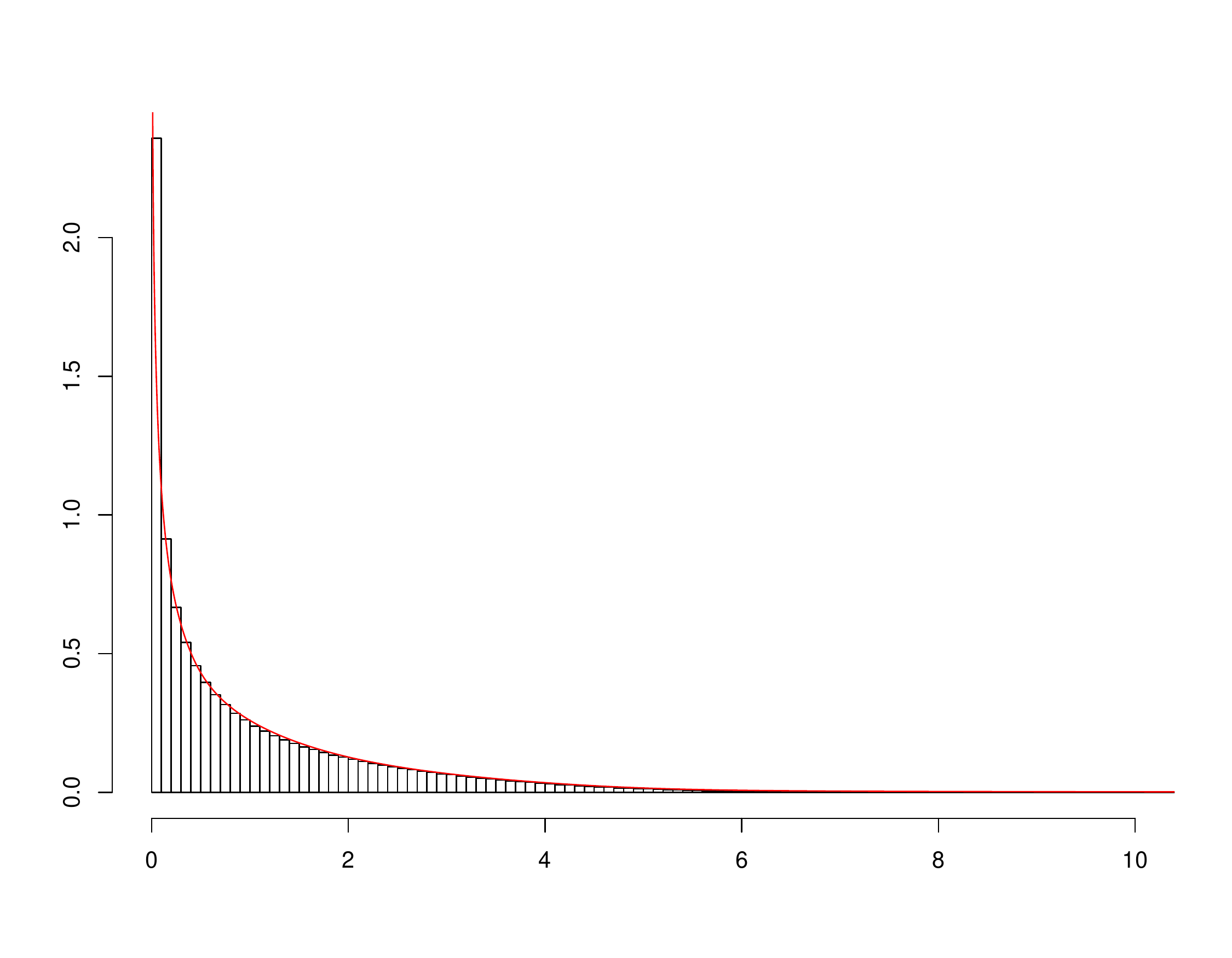}
         \caption{Comparison between the theoretical value of the density $\upsilon\circ(x^2)^{-1}(x)$ and the empirical histogram computed through a sample of simulated empirical covariance matrices $R_N$
         as defined in the introduction. For both plots, $q=1$ but stock prices follow multifractal random walks with intermittency parameter $\gamma^2=1/4$ in the upward figure, $\gamma^2=1/2$ in the downward figure. }\label{intermittence}
\end{figure}

The results are as follows. In Fig. \ref{intermittence}, we show the comparison between 
the theoretical value of the density $\upsilon\circ(x^2)^{-1}(x)$ (computed numerically as described above) 
and an empirical histogram of the eigenvalues of a sample of simulated covariance matrices $R_N$ (defined in the introduction) for $N=1024$ and $q=1$.  The upward plot is done with an intermittency parameter $\gamma^2 = 1/4$ and an integral scale $\tau=1/4$. The agreement is excellent as expected from Theorems \ref{main2},  \ref{main1} and  \ref{main3}. The downward figure is done for an intermittency parameter $\gamma^2 = 1/2$ 
and an integral scale $\tau=1/4$,
suggesting that our prediction remains true for $\gamma^2>1/3$ (see conjecture \ref{conj} which also covers the case $\gamma^2 \in [1,2[$). 

\begin{figure}[h!btp] 
	\center
         \includegraphics[scale=0.6]{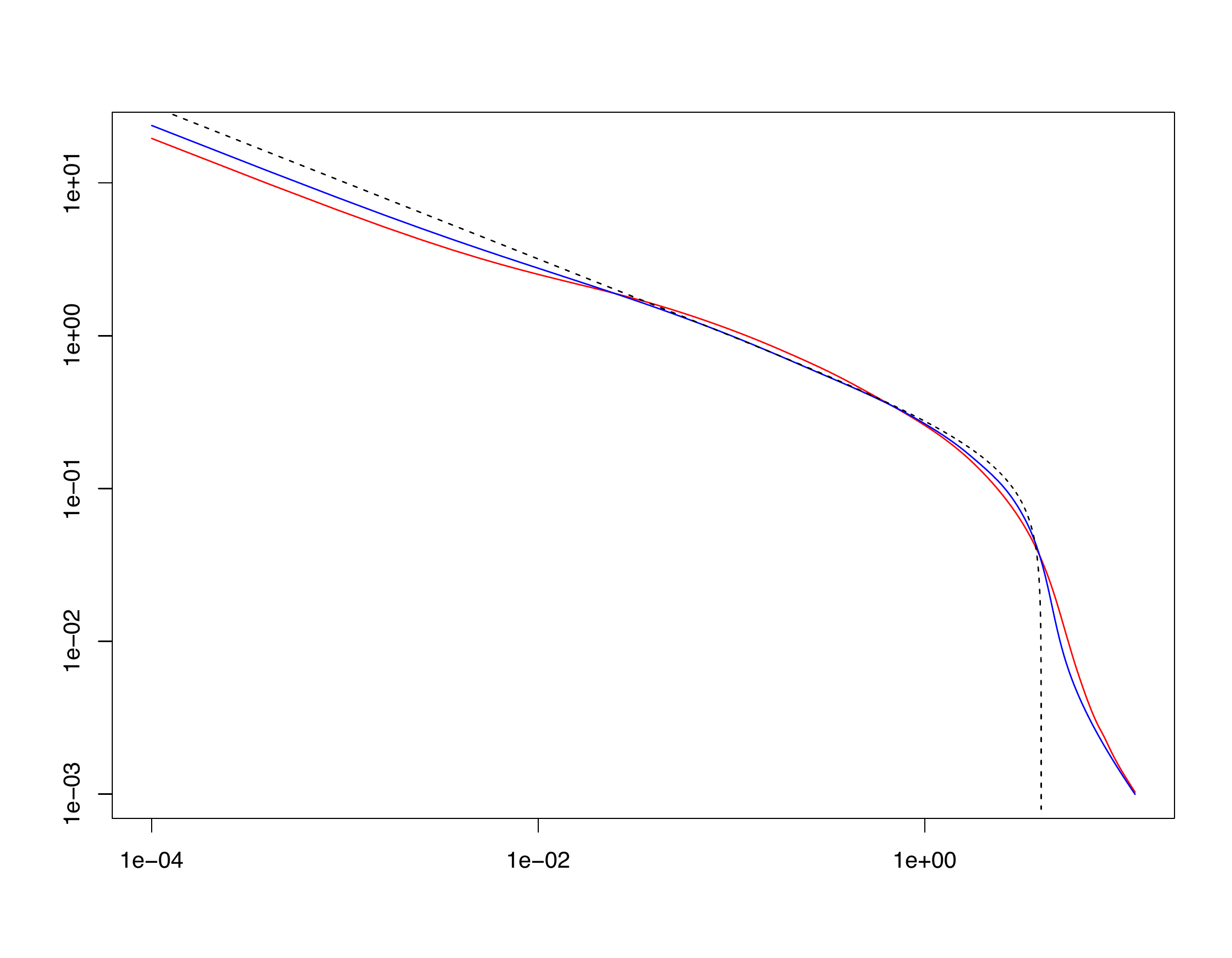}
         \caption{Log-log plot of the density $\upsilon\circ(x^2)^{-1}$ with $q=1$, $\tau = 1/4$ for three different intermittency parameter: $\gamma^2 = 0$ (black dashed line), $\gamma^2=1/4$ (blue line) and 
         $\gamma^2=1/2$ (red line). }\label{comp_courbe}
\end{figure}

In Fig. \ref{comp_courbe}, we represent three curves (axis are in log-log) corresponding to the theoretical density  $\upsilon\circ(x^2)^{-1}(x)$ for a parameter $q=1$, an integral scale $\tau=1/4$ 
and for three different values of $\gamma^2$. The black dashed curve 
corresponds to $\gamma^2=0$, which in fact is the Marcenko-Pastur case: asset prices are following independent Brownian motions with a trivial constant volatility process. In this case, the support is compact and the 
right edge of the spectrum is known to be equal to $4$. The blue curve corresponds to an intermittency parameter equal to $1/4$ and the red curve is for $\gamma^2=1/2$.
In this way, we see precisely the distortion of the spectrum induced by the auto-correlated volatility process. The most interesting part for applications is certainly about the tails of the distribution: the higher the intermittency 
parameter $\gamma^2$ is, the heavier the tail of the distribution is.    

\begin{figure}[h!btp] 
	\center
         \includegraphics[scale=0.6]{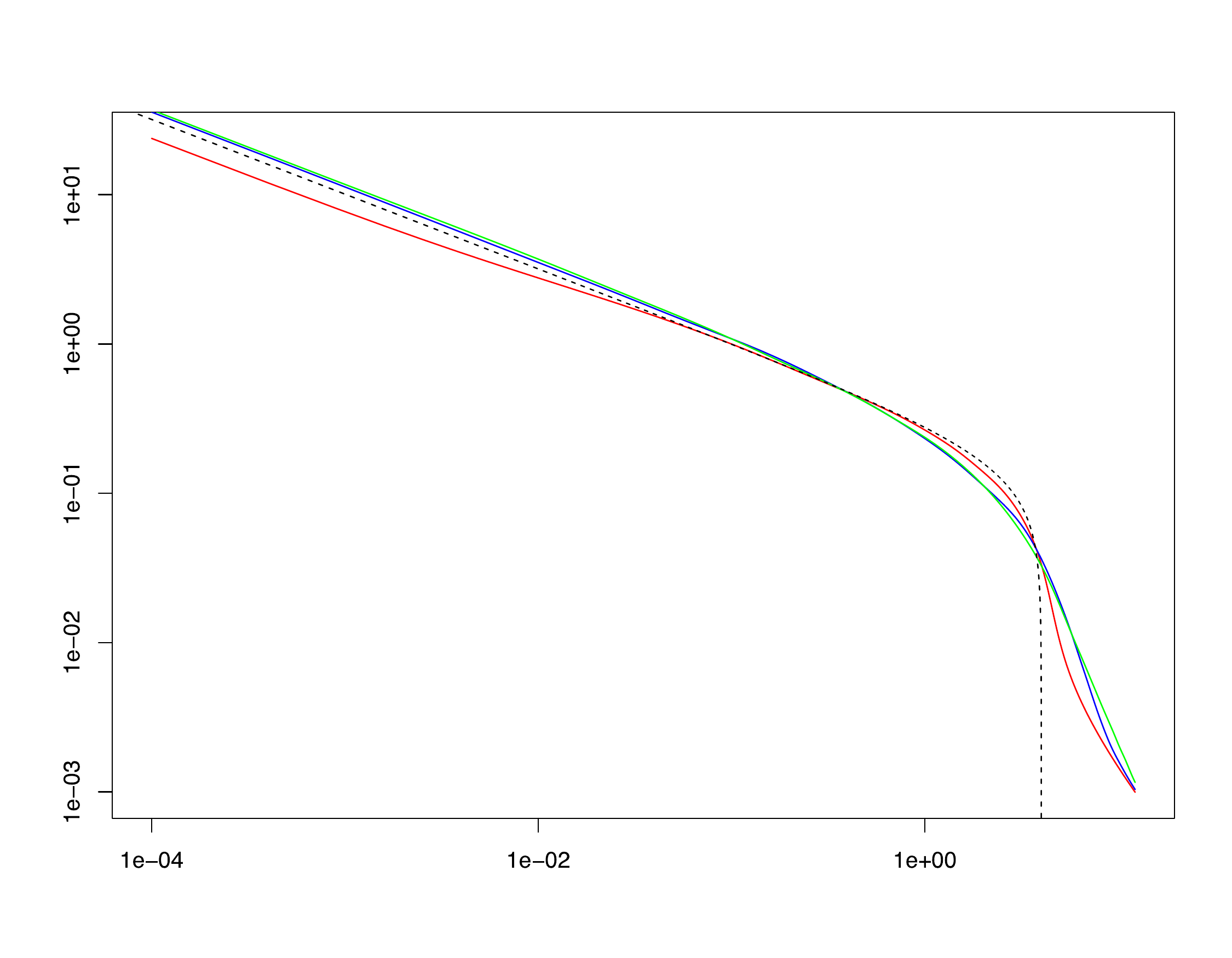}
         \caption{Log-log plot of the density $\upsilon\circ(x^2)^{-1}$ with $q=1,\gamma^2=1/4$ for four different integral scales $\tau$: $\tau = 0$ (black dashed line), $\tau=1/4$ (red line), 
         $\tau=1$ (blue line) and $\tau=2$ (green line).}\label{comp_courbe_tau}
\end{figure}

In Fig. \ref{comp_courbe_tau}, we represent four curves corresponding to the thoeretical density $\upsilon\circ(x^2)^{-1}(x)$ but varying the integral scale $\tau$ instead of the intermittency parameter $\gamma^2$. 
We chose for this plot $q=1$ and $\gamma^2=1/4$ and represented the density $\upsilon\circ(x^2)^{-1}(x)$ for $\tau=0$ (corresponding to the trivial MP case) and for $\tau=1/4,1,2$. The result on the right tail 
of the distribution is the following: the higher the integral scale is, the heavier the right tail of the distribution is. As mentionned above, large integral scale corresponds to 
measuring price variations on small scales. On small scales, it is known that price variations will have distribution with larger kurtosis than price variations on larger scales and therefore it was expected to find 
heavier right tail distribution for the spectral distribution of the corresponding covariance matrix.    

\section{Proofs of the main results}\label{proofs}

In this section, we give the proofs of theorems \ref{main2}, \ref{main1} and \ref{main3}.
The proof of Theorem \ref{main4} is very similar and we will not explain it in every detail, except for the final part where we establish the second equation of the system 
in Theorem \ref{main4} verified by $K_z$. We will give the details for this part of the proof in the appendix. 
The proof of theorem \ref{main3lognorm} is an easy adaptation of our proofs for theorems \ref{main2}, \ref{main1} and \ref{main3}; it is left to the reader. 
Furthermore, the proofs are very similar when $q=1$ or when $q<1$. For the sake of clarity, we assume $T=N$ and hence $q=1$ in the proofs that follow.

Hence, in the following, we will suppose (unless otherwise stated) that:
\begin{align*}
\varphi(q)& = -iq\frac{\gamma^2}{2} + q^2 \frac{\gamma^2}{2}, \\
\psi(q)& = \varphi(-iq), \\
\zeta(q)& = (1+\frac{\gamma^2}{2})q + q^2 \frac{\gamma^2}{2}, \\
\gamma^2& < \frac{1}{3},
\end{align*}
and $M$ will be the MRM whose structure exponent is $\zeta$ (see section \ref{cons:loglevy} for a reminder).

Our approach to show the convergence of $\E[L^{1,z}_N]$ and $\E[L^{2,z}_N]$ consists in proving tightness and characterizing uniquely
the possible limit points. The
classical Schur complement formula is our basic linear algebraic tool to study $\E[L^{1,z}_N]$ and $\E[L^{2,z}_N]$ recursively on
the dimension $N$, as is usual when the resolvent method is used. The original part of our proof is that we apply the Schur complement formula two times in a row to find the second equation
of the system in theorem \ref{main3} involving the limit point $K_z(x)$ of the measure $\E[L^{1,z}_N]$. We will also show that the limit points of the two complex measures $\E[L^{1,z}_N]$ and $\E[L^{2,z}_N]$
satisfy a fixed point system (written in theorem \ref{main3}). 

We begin by showing tightness.

\subsection{Tightness of the complex measures $\E[L^{1,z}_N], \E[L^{2,z}_N]$ and limit points}
\begin{lemma}
The two families of complex measures $(\E[L^{i,z}_N])_{N \in \N}, i=1,2$ are tight and bounded in total variation.  
\end{lemma}
 
\proof Let us present the proof for $(\E[L^{1,z}_N])_{N \in \N}$; the other proof is similar. 

One has, for each $N \in \N$:
\begin{equation}
\mid \E[L^{1,z}_N] \mid [0,1] = \frac{1}{N} \sum_{k=1}^{N} \mid \E[G_{N}(z)_{kk}] \mid \leq \frac{1}{\mid \Im(z) \mid},
\end{equation}
and so the family of complex measures $(\E[L^{1,z}_N])_{N \in \N}$ is bounded in total variation. 
It is obviously tight since the support of all the complex measures in the family is included in $[0,1]$, which is a compact set.
\hfill$\square$

Using Prokhorov's theorem, we know that those two families of complex measures are sequentially compact
in the space of complex Borel measure on $[0,1]$ equipped with the topology of weak convergence.
In particular, there exists a subsequence such that, for all bounded continuous function $f$, one has, when $N$ goes to $+\infty$ along this subsequence:
\begin{equation}
\E\left[ L^{1,z}_N(f) \right] \to \int_0^1 f(x) \mu_z^1(dx).  
\end{equation}

\begin{lemma}
The complex measure $\mu_z^1(dx)$ has Lebesgue density; more precisely, there exists a bounded measurable function $K_z(x)$ such that:
\be
\mu_z^1(dx) = K_z(x) dx.
\ee
\end{lemma}

\proof
One has: 
\begin{align}
\left| \E\left[  L^{1,z}_N(f) \right] \right| &\leq \frac1N \sum_{k=1}^N |f(k/N)| \E\left[G_{N}(z)_{kk}\right] \\
                                           &\leq \frac{1}{|\Im(z)|} \frac1N \sum_{k=1}^N |f(k/N)|
\end{align}

Letting $N \to +\infty$ along a subsequence, one obtains:
\be
\left| \int_0^1 f(x) \mu_z^1(dx) \right| \leq \frac{1}{|\Im(z)|} \int_0^1 |f(x)| dx.
\ee
This proves the lemma. 
\qed

Thus, there exists a subsequence such that, as $N$ tends to $+\infty$ along this subsequence:
\begin{equation}
\E\left[ L^{1,z}_N(f) \right] \to \int_0^1 f(x) K_z(x) dx.  
\end{equation}

\begin{lemma}
There exists a subsequence and a constant $\mu_z^2 \in \C$ such that, as $N$ goes to $+\infty$ along this subsequence:
\be
\E\left[ L^{2,z}_N(f) \right] \to  \mu_z^2 \int_0^1 f(x) dx.
\ee
\end{lemma}

\proof It is easy to see that the $G_{N}(z)_{kk}, k=N+1,\dots,N$ are identically distributed. In particular, these variables have the same mean $\mu_z^2(N)$.
One has, for all $N$:
\be
|\mu_z^2(N)| \leq \frac{1}{|\Im(z)|}. 
\ee
So there exists a subsequence and a complex number $\mu_z^2$ such that, as $N$ goes to $+\infty$ along this subsequence, $\mu_z^2(N) \to \mu_z^2$. 
One thus obtains, as $N$ goes to $+\infty$ along this subsequence: 
\be
\E\left[ L^{2,z}_N(f) \right] \to  \mu_z^2 \int_0^1 f(x) dx.
\ee
\hfill$\square$

Following the classical method as in \cite{BenArousGuionnet}, \cite{GuionnetDemboBelinshi}, \cite{BandMatricesVariance}, we will show in the following that the limit point $\mu_z^2$ and $K_z(x)$ are defined uniquely 
and do not depend on the subsequence. We will first recall some preliminary results on resolvents.  

\subsection{Preliminary results on resolvents}
We first recall the following standard and general result; the next lemmas of this section are also standard but are applied to our particular case.
\begin{lemma}\label{lem:im}
Let $A$ be a symmetric real valued matrix of size $N$. For $z\in\C\setminus \R$, let us denote by $G(z)$ the matrix
\be
G(z)=(z-A)^{-1}.
\ee 
For $z\in\C\setminus \R$ and $k\in \{1,\dots,N\}$, we have
\be\label{imz}
\Im(z)\Im(G(z)_{kk})< 0\quad \text{ and }|G(z)_{kk}|\leq \frac{1}{|\Im(z)|}. 
\ee
In particular, if $F\subset\{1,\dots,N\}$ is a finite set and $(a_i)_{i\in F}$ a finite sequence of positive number, then: 
\be\label{imz2}
\frac{\Im\Big(z-\sum_{i\in F}a_iG(z)_{ii}\Big)}{\Im(z)}\geq 1.
\ee
and we also have:
\be
\frac{1}{\left|z-\sum_{i\in F}a_iG(z)_{ii}\right|}\leq \frac{1}{|\Im(z)|}.
\ee
\end{lemma}

\proof Write $A=\bar{U}^tDU$ where $D$ is a diagonal matrix with diagonal real entries $(\lambda_i)_{1\leq i\leq N}$. Then 
$$G(z)_{kk}=\sum_{i=1}^N|U_{ki}|^2\frac{1}{z-\lambda_i}.$$
Since $\Re\Big(\frac{1}{z-\lambda_i}\Big)=\frac{\Re(z)-\lambda_i}{(\Re(z)-\lambda_i)^2+\Im(z)^2}$ and $\Im\Big(\frac{1}{z-\lambda_i}\Big)=\frac{-\Im(z)}{(\Re(z)-\lambda_i)^2+\Im(z)^2}$ the relation \eqref{imz} follows. 
It is then straightforward to derive \eqref{imz2} from \eqref{imz}. \qed

For $i=1, \dots, N$, let $X_N^{(i)}=\left(X_N(kl)\right)_{k,l\not = i}$ be the matrix obtained from $X_N$ by taking off the $i$-th column and row. 
Define, also for $i=1, \dots, 2N $ the $(2N-1)\times (2N-1)$ matrix $A^{(i)}_{N}(z)$ obtained from $A_{N}(z)$ by taking off the $i$-th column and row.
In particular, for $i=1, \dots, N$, 
\begin{equation*}
 A^{(N+i)}_{N}(z)=
 \begin{pmatrix}
         z I_{N} & - \! {}^{t} \! X_N^{(i)}\\
         -X_N^{(i)} & z I_{N-1}
\end{pmatrix},  
\end{equation*}
For $i=1, \dots, 2N $, set: 
\be
G^{(i)}_{N}(z) = (A^{(i)}_{N}(z))^{-1}.
\ee
Let now $\hat{X}_N^{(i)}$ denote the matrix $X_N$ with the $i$-th column and row set to $0$ and $\hat{A}^{(i)}_{N}(z)$ denote the matrix $A_{N}(z)$ with the $i$-th column and row set to $0$ excepted the diagonal term.
Again we have, for $i=1, \dots, N$:
\begin{equation*}
\hat{A}^{(N+i)}_{N}(z)=
 \begin{pmatrix}
         z I_{N} & -{}^{t} \! \hat{X}_N^{(i)} \\
         -\hat{X}_N^{(i)} & z I_{N}
\end{pmatrix},  
\end{equation*}
For $i=1, \dots, 2N$, set: 
\be
\hat{G}^{(i)}_{N}(z) = (\hat{A}^{(i)}_{N}(z))^{-1}.
\ee
In the paper, we will also use the terms $A^{(k,i)}_{N}(z),G^{(k,i)}_{N}(z),\hat{A}^{(k,i)}_{N}(z),\hat{A}^{(k,i)}_{N}(z) $. The double superscript just means that you make the operations described above to 
the rows and columns $i$ and $k$. 

\begin{lemma}\label{termDown}
For all $k \in \{ 1, \dots, N\}$ and all $t \not = N+k$, one has:
\be
\E\left[\left| G_{N}(z)_{tt} - \hat{G}^{(N+k)}_{N}(z)_{tt} \right| \right] \leq \frac{1}{\sqrt{N} |\Im(z)|^2}.
\ee
\end{lemma}

\proof
Multiply the identity:
\be
 \hat{A}^{(N+k)}_{N}(z) - A_{N}(z)= \hat{A}^{(N+k)}_{N}(0) - A_{N}(0)
\ee
to the left by $G_{N}(z)$ and to the right by $\hat{G}^{(N+k)}_{N}(z)$ to obtain
\be
G_{N}(z) - \hat{G}^{(N+k)}_{N}(z) = G_{N}(z) ( \hat{A}^{(N+k)}_{N}(0) - A_{N}(0) ) \hat{G}^{(N+k)}_{N}(z).
\ee
Then one has:
\begin{align}
G_{N}(z)_{tt} - \hat{G}^{(N+k)}_{N}(z)_{tt} &= \left(G_{N}(z) ( \hat{A}^{(N+k)}_{N}(0) - A_{N}(0) ) \hat{G}^{(N+k)}_{N}(z)\right)_{tt} \\
                                                  &= \hat{G}^{(N+k)}_{N}(z)_{N+k,t} \sum_{i=1}^{N} G_{N}(z)_{ti} r_k(i) \\
                                                  &+ G_{N}(z)_{t,N+k} \sum_{j=1}^{N} r_k(j) \hat{G}^{(N+k)}_{N}(z)_{jt} \\
                                                  &= G_{N}(z)_{t,N+k} \sum_{j=1}^{N} r_k(j) \hat{G}^{(N+k)}_{N}(z)_{jt}
\end{align}
where we have noticed that, for all $t \not = N+k, \hat{G}^{(N+k)}_{N}(z)_{N+k,t} = 0$.

Therefore, we find that:
\be
\E\left[ \left| G_{N}(z)_{tt} - \hat{G}^{(N+k)}_{N}(z)_{tt} \right| \right] \leq \E\left[ \left| G_{N}(z)_{t,N+k} \right|^2 \right]^{1/2} 
\E\left[ \left| \sum_{j=1}^{N} r_k(j) \hat{G}^{(N+k)}_{N}(z)_{jt} \right|^2 \right]^{1/2} 
\ee
by Cauchy-Schwartz's inequality. Using then the independence of $r_k(j)$ and $\hat{G}^{(N+k)}_{N}(z)$, we get:
\begin{align*}
\E\left[ \left| G_{N}(z)_{tt} - \hat{G}^{(N+k)}_{N}(z)_{tt} \right| \right] &\leq \E\left[ \left| G_{N}(z)_{t,N+k} \right|^2 \right]^{1/2} 
\E\left[ r_k(1)^2 \right]^{1/2} \E\left[  \sum_{j=1}^{N} \left| \hat{G}^{(N+k)}_{N}(z)_{jt} \right|^2 \right]^{1/2} \\
& \leq \frac{1}{\sqrt{N} |\Im(z)|^2}.
\end{align*}
The proof is complete.\qed

\begin{lemma}\label{termUP}
There exists a constant $C>0$ such that, for all $k \in \{ 1, \dots, N\}$ and all $t \not = k$:
\be
\E\left[\left| G_{N}(z)_{tt} - \hat{G}^{(k)}_{N}(z)_{tt} \right| \right] \leq \frac{C}{|\Im(z)|^2} \frac{1}{N^{\frac{1-\gamma^2}{4}}} .
\ee
\end{lemma}

\proof
Again, we start from the relation:
$$
G_{N}(z) - \hat{G}^{(k)}_{N}(z) = G_{N}(z) ( \hat{A}^{(k)}_{N}(0) - A_{N}(0) ) \hat{G}^{(k)}_{N}(z).
$$
Thus we have
\begin{align}
G_{N}(z)_{tt} - \hat{G}^{(k)}_{N}(z)_{tt} &= \left(G_{N}(z) ( \hat{A}^{(k)}_{N}(0) - A_{N}(0) ) \hat{G}^{(k)}_{N}(z)\right)_{tt} \\
                                                  &= \hat{G}^{(k)}_{N}(z)_{k,t} \sum_{i=N+1}^{N} G_{N}(z)_{ti} r_i(k) \\
                                                  &+ G_{N}(z)_{t,k} \sum_{j=1}^{N+1} r_j(k) \hat{G}^{(k)}_{N}(z)_{jt} \\
                                                  &= G_{N}(z)_{t,k} \sum_{j=1}^{N+1} r_j(k) \hat{G}^{(k)}_{N}(z)_{jt}
\end{align}
where we have noticed that, for all $t \not = k, \hat{G}^{(k)}_{N}(z)_{k,t} = 0$.

Therefore, we find that:
\be
\E\left[ \left| G_{N}(z)_{tt} - \hat{G}^{(k)}_{N}(z)_{tt} \right| \right] \leq \E\left[ \left| G_{N}(z)_{t,k} \right|^2 \right]^{1/2} 
\E\left[ \left| \sum_{j=1}^{N} r_j(k) \hat{G}^{(k)}_{N}(z)_{jt} \right|^2 \right]^{1/2} 
\ee
by Cauchy-Schwartz's inequality. We want to expand the square in the above expression. To that purpose, we first observe that, conditionally to the $M^i$, the variables $(r_j(k))_j$ are independent from $\hat{G}^{(k)}_{N}(z) $ and centered. Hence we have for $j\not = j'$,
$$\E\left[   r_j(k)r_{j'}(k) \hat{G}^{(k)}_{N}(z)_{jt} \hat{G}^{(k)}_{N}(z)_{j't}   \right] =0.$$
Thus we get:
\begin{align*}
\E\left[ \left| G_{N}(z)_{tt} - \hat{G}^{(k)}_{N}(z)_{tt} \right| \right] &\leq \E\left[ \left| G_{N}(z)_{t,k} \right|^2 \right]^{1/2} 
\left( \sum_{j=1}^{N+1} \E\left[ r_j(k)^2 \left|\hat{G}^{(k)}_{N}(z)_{jt} \right|^2 \right]\right)^{1/2}  \\
&\leq \E\left[ \left| G_{N}(z)_{t,k} \right|^2 \right]^{1/2} \left( \sum_{j=1}^{N+1} \E\left[ r_j(k)^4\right]^{1/2} \E\left[\left|\hat{G}^{(k)}_{N}(z)_{jt} \right|^4 \right]^{1/2} \right)^{1/2} \\
& \leq \frac{\E[r_1(k)^4]^{1/4}}{|\Im(z)|} \left(\sum_{j=1}^{N+1} \E\left[\left|\hat{G}^{(k)}_{N}(z)_{jt} \right|^4 \right]^{1/2} \right)^{1/2} \\
& \leq \frac{\E[r_1(k)^4]^{1/4}}{|\Im(z)|}(N+1)^{1/4} \left(   \sum_{j=1}^{N+1} \E\left[\left|\hat{G}^{(k)}_{N}(z)_{jt} \right|^4 \right] \right)^{1/4} 
\end{align*}
Now we use the scaling properties of the MRM to obtain, for some positive constant $C$,
$$\E[r_j(k)^4]=3\E[M(0,\frac{1}{N})^2]\leq CN^{-\zeta(2)}.$$
Furthermore, by using Lemma \ref{SommeLigne} which assures that, almost surely:
\be
\sum_{j=1}^{N+1} \left|\hat{G}^{(k)}_{N}(z)_{jt} \right|^2  \leq \frac{1}{|\Im(z)|^2}
\ee
and the fact that:
\be
\sum_{j=1}^{N+1} \left|\hat{G}^{(k)}_{N}(z)_{jt} \right|^4  \leq \left( \sum_{j=1}^{N+1} \left|\hat{G}^{(k)}_{N}(z)_{jt} \right|^2 \right)^2,
\ee
we finally obtain
\begin{align*}
\E\left[ \left| G_{N}(z)_{tt} - \hat{G}^{(k)}_{N}(z)_{tt} \right| \right]   
&\leq \frac{C}{|\Im(z)|^2} \big(\frac{1}{N}\big)^{\frac{\zeta(2)-1}{4}}.
\end{align*}
It just remains to check that $\zeta(2) = 2 - \gamma^2$. 
\qed

\begin{lemma}\label{Gchap}
For each $k\in\{1,\dots,2N\}$, if $t \not = k$, then 
\be 
G^{(k)}_{N}(z)_{tt} = \hat{G}^{(k)}_{N}(z)_{tt}, 
\ee 
and if $t = k$, then $\hat{G}^{(k)}_{N}(z)_{k,k} = z^{-1}$. 
\end{lemma}

\proof
It is straightforward to see that the two matrices $G^{(k)}_{N}(z)$ and $\hat{G}^{(k)}_{N}(z)$ have the same eigenvalues except that $\hat{G}^{(k)}_{N}(z)$ has one more zero eigenvalue. 
In addition, the eigenvectors look also very similar since you can obtain $2N$ eigenvectors of $\hat{G}^{(k)}_{N}(z)$ by adding a zero entry to the eigenvectors of $G^{(k)}_{N}(z)$ (between the entries
$k-1$ and $k$). The last eigenvector of $\hat{G}^{(k)}_{N}(z)$ is the vector of $\R^{N}$ for which all entries are zero except the entry number $k$. 

Now observe that with $G^{(k)}_{N}(z) = U \text{diag}(z-\lambda) U^*$ and $\hat{G}^{(k)}_{N}(z) = V \text{diag}(z-\tilde{\lambda}) V^*$, 
\begin{align}
G^{(k)}_{N}(z)_{tt} &= \sum_{i=1}^{2N} |u_{ti}|^2 \frac{1}{z-\lambda_i} \\
\hat{G}^{(k)}_{N}(z)_{tt} &= \sum_{i=1}^{N} |v_{ti}|^2 \frac{1}{z-\tilde{\lambda_i}}.                        
\end{align}

The result follows since, for $t \not = k$,
\be
\sum_{i=1}^{2N-1} |u_{ti}|^2 \frac{1}{z-\lambda_i} = \sum_{i=1}^{2N } |v_{ti}|^2 \frac{1}{z-\tilde{\lambda_i}}
\ee
and, for $t = k, \hat{G}^{(k)}_{N}(z)_{k,k} = z^{-1}.$
\qed

\begin{lemma}\label{imK}
For all   $z\in\C$ and Lebesgue almost every point $x\in [0,1]$, we have
\begin{equation}
\Im(z)\Im(K_z(x))\leq 0
\end{equation}
and 
\begin{equation}
|\Im(K_z(x))| \leq \frac{1}{\Im(z)}
\end{equation}
\end{lemma}

\proof  This is a straightforward consequence of Lemma \ref{lem:im}. Indeed, we have for all positive continuous function $f$ on $[0,1]$ and $N\in\N$:
$$\Im(z)\Im\Big(\int_0^1f(x)\E[L^{1,z}_N](dx)\Big)\leq 0.$$
We pass to the limit as $N$ goes to $\infty$ along some suitable subsequence and obtain:
$$\Im(z)\Im\Big(\int_0^1f(x)K_z(x)\,dx\Big)\leq 0.$$
The result follows.\qed

\subsection{Concentration inequalities}
This lemma is adapted to our case from Lemma 5.4 in \cite{BenArousGuionnet}. 
\begin{lemma} \label{concentration}
Let $f:[0,1] \to \R$ be a bounded measurable function. 
For each $i \in \{1,2\}$, we have the following concentration results:
\begin{equation}
\E\left[ \mid L^{i,z}_N(f) - \E[L^{i,z}_N(f)]  \mid^2 \right] \leq  \frac{8}{N} \frac{\mid\mid f \mid \mid^2_{\infty}}{\mid \Im{z} \mid^4}.
\end{equation}
\end{lemma}

\proof  Define two functions $F^1_N$ and $F^2_N$ such that:
\begin{align}
F^1_N\left(\left(X_{ij}^{(N)}\right)_{1\leq j \leq N+1}, 1\leq i \leq N\right) &= \frac{1}{N} \sum_{k=1}^{N} f\left(\frac{k}{N}\right) G_{N}(z)_{kk}  \\
F^2_N\left(\left(X_{ij}^{(N)}\right)_{1\leq j \leq N+1}, 1\leq i \leq N\right) &= \frac{1}{N} \sum_{k=1}^{N+1}  f\left(\frac{k}{N+1}\right) G_{N}(z)_{k+N,k+N}
\end{align}

We will prove the Lemma for $L^{1,z}_N$; the proof for $L^{2,z}_N$ is a straightforward adaptation. 

Let, for $k \in \{ 1, \dots, N+1 \}$, 
\begin{equation}
\mathcal{F}_k = \sigma\left(\left(X_{ij}^{(N)}\right)_{1\leq j \leq N}, 1\leq i \leq k\right)
\end{equation}

If $P$ denotes the law of the vector $\left(X_{1j}^{(N)}\right)_{1 \leq j \leq N}$,
\begin{align*}
&\E\left[ \mid F^1_N - \E[F^1_N]  \mid^2 \right] \\
&= \sum_{i=0}^{N} \E\left[ \mid \E[ F^1_N\mid \mathcal{F}_{i+1}] - \E[F^1_N\mid \mathcal{F}_i]  \mid^2 \right] \\
&= \sum_{i=0}^{N} \int \mid \int \left( F_N(x_1, x_2,\dots, x_{i+1}, y_{i+2}, \dots, y_{N+1}) - F_N(x_1, x_2,\dots, x_{i}, y_{i+1}, \dots, y_{N+1}) \right) dP^{\otimes {N+1}}(y) \mid^2 \\ 
&dP^{\otimes {i+1}}(x) \\
&\leq  \sum_{i=0}^{N} \int \mid \int \left( F_N(x_1, x_2,\dots,x_i, x_{i+1}, x_{i+2}, \dots, x_{N+1}) - F_N(x_1, x_2,\dots, x_{i}, y, x_{i+2}, \dots, x_{N+1}) \right) dP(y) \mid^2 \\ 
&dP^{\otimes {N+1}}(x) \\
& \leq \sum_{i=0}^{N}\sup_{\R^{(N+1)^2}} \mid\mid \nabla_{x_{i+1}} F_N \mid\mid^2 \int \mid\mid x-y \mid\mid^2 dP^{\otimes 2}(x,y). 
\end{align*}

The quantity $\nabla_{x_{i+1}} F^1_N$ refers to the gradient of $F^1_N$ in the direction of the vector $x_{i+1}$.

If we consider a couple of processes $(\widetilde{B}^1,\widetilde{M}^1)$ independent from $(B^1,M^1)$ with the same law,  it is easy to see that:
\begin{align*}
 \int \mid\mid x-y \mid\mid^2 dP \otimes dP(x,y) & = \sum_{j=1}^N\E\Big[(B^1_{M^1(0,\frac{j}{N})}-B^1_{M^1(0,\frac{j-1}{N})}-\widetilde{B}^1_{\widetilde{M}^1(0,\frac{j}{N})}+\widetilde{B}^1_{\widetilde{M}^1(0,\frac{j-1}{N})})^2\Big].\\
 &=2-2\sum_{j=1}^N\E\Big[(B^1_{M^1(0,\frac{j}{N})}-B^1_{M^1(0,\frac{j-1}{N})})(\widetilde{B}^1_{\widetilde{M}^1(0,\frac{j}{N})}\widetilde{B}^1_{\widetilde{M}^1(0,\frac{j-1}{N})})\Big]\\
 &= 2.
\end{align*}
In our case, we have, for $i \in \{ 1, \dots, N+1 \}, j \in \{ 1, \dots, N \}$:
\begin{equation}
\frac{\partial G_{N}(z)_{kk}}{\partial X_{ij}} = G_{N}(z)_{k,j} G_{N}(z)_{N+i,k} + G_{N}(z)_{k,N+i} G_{N}(z)_{j,k}
\end{equation}

Thus, 
\begin{equation}
\nabla_{x_{i+1}} F_N = \frac{1}{N} \sum_{k=1}^{N} f\left(\frac{k}{N}\right) \nabla_{x_{i+1}} G_{N}(z)_{kk}
\end{equation}
It is now plain to compute:
\begin{align*}
\mid\mid \nabla_{x_{i+1}} F_N \mid\mid^2 &= \frac{1}{N^2} \sum_{j=1}^{N} \mid \left( G_{N}(z) D^1(f) G_{N}(z) \right)_{N+i+1,j} + \left( G_{N}(z) D^1(f) G_{N}(z) \right)_{j,N+i+1} \mid^2 
\end{align*}
where $D^1(f)$ is the $(2N)$-dimensional diagonal matrix of entries:
\begin{align*}
D^1(f)_{kk} = f\left(\frac{k}{N} \right) 1_{\{1\leq k \leq N\}}.
\end{align*}
 
One thus has:
\begin{align*}
\mid\mid \nabla_{x_{i+1}} F^1_N \mid\mid^2 &=\frac{4}{N^2}  \sum_{j=1}^{N} \mid \left( G_{N}(z) D^1(f) G_{N}(z) \right)_{N+i+1,j} \mid^2 \\
& \leq \frac{4}{N^2}  \sum_{j=1}^{2N} \mid \left( G_{N}(z) D^1(f) G_{N}(z) \right)_{N+i+1,j} \mid^2 \\
&\leq \frac{4}{N^2} \frac{\mid\mid f \mid\mid^2_\infty }{\mid \Im{z} \mid^4}.
\end{align*}
where, in the last line, we used lemma \ref{SommeLigne} and the fact that the matrix $G_{N}(z) D^1(f) G_{N}(z)$ has a spectral radius smaller than $\mid\mid f \mid\mid_\infty /\mid \Im{z} \mid^2$.

Finally,
\begin{equation}
\E\left[ \mid F^1_N - \E[F^1_N]  \mid^2 \right] \leq  \frac{8}{N} \frac{\mid\mid f \mid\mid^2_\infty}{\mid \Im{z} \mid^4}.
\end{equation}
\hfill$\square$

We also prove the following lemma:
\begin{lemma}\label{concsup}
For all $\alpha>1$ such that $\zeta(2\alpha)>1$, we have
\begin{align}
\E\left[\left|  \sum_{t=1}^{N}r_k(t) ^2  \left(\hat{G}^{(N+k)}_{N}(z)_{tt} - \E[\hat{G}^{(N+k)}_{N}(z)_{tt}] \right) \right|  \right] 
&\leq \frac{C (\ln N)^2}{N^{\frac{\zeta(2 \alpha)-1}{\alpha}} |\Im(z)|^4}
\end{align}
for some positive constant $C$ independent from $N,z,k$.
\end{lemma}

\proof Notice that $(r_k(t))_t$ and $\hat{G}^{(N+k)}_{N}(z)$ are independent. Hence, by conditioning with respect to the process $(r_k(t))_t$, we can argue along the same lines as in the previous lemma with $r_k(t)$ 
instead of $\frac{1}{N}f(\frac{t}{N})$ and we get the formula:
\begin{align*}
 \E\left[\left|  \sum_{t=1}^{N}r_k(t) ^2  \left(\hat{G}^{(N+k)}_{N}(z)_{tt} - \E[\hat{G}^{(N+k)}_{N}(z)_{tt}] \right) \right|^2 \right] 
&\leq    \frac{8}{|\Im(z)|^4}\E[\sup_tr_k(t)^4].
\end{align*}
We conclude with Proposition \ref{sup} in the appendix .\qed

In the following, we fix $\alpha >1$ such that $\zeta(2\alpha)>1$ (because of the expression of $\zeta$ and the inequality $\gamma^2 < 1/3$, it is clear that such a number exists). 


\subsection{The system verified by the limit point $\mu_z^2$ and $K_z(x)$: first equation}\label{section:firstequation}
From the Schur complement formula (see e.g. Lemma 4.2 in \cite{BenArousGuionnet} for a reminder), one has for $k \in \{1, \dots, N\}$:
\be\label{defgnk}
G_{N}(z)_{N+k,N+k} = \left[ z - \sum_{s,t=1}^{N} r_k(s) r_k(t) G^{(N+k)}_{N}(z)_{st} \right]^{-1}
\ee
Using Lemma \ref{lem:diag1}, one can write:
\be
G_{N}(z)_{N+k,N+k} = \left[ z - \sum_{t=1}^{N} r_k(t)^2 G^{(N+k)}_{N}(z)_{tt} + \epsilon^1_{N,k}(z) \right]^{-1}
\ee
where $\epsilon^1_{N,k}(z)$ is a complex valued random variable for which there exists $C>0$ such that for all $N \in \N$ and $1\leq k \leq N$, 
\be
\E[|\epsilon^1_{N,k}(z)|^2] <  \frac{C}{N^{1-\gamma^2}}.
\ee

By using Lemma \ref{Gchap}, we can write:
\begin{align}
G_{N}(z)_{N+k,N+k} &= \left[ z - \sum_{t=1}^{N} r_k(t)^2 \hat{G}^{(N+k)}_{N}(z)_{tt} + \epsilon^{1}_{N,k}(z) \right]^{-1} .
\end{align}

Lemma \ref{concsup} applied to $\alpha>1$ such that $\zeta(2\alpha)>1$ yields:
\begin{align}
&\E\left[\left|  \sum_{t=1}^{N}r_k(t)^2   \left(\hat{G}^{(N+k)}_{N}(z)_{tt} - \E[\hat{G}^{(N+k)}_{N}(z)_{tt}] \right) \right|^2 \right] \leq \frac{C (\ln N)^2}{N^{\frac{\zeta(2\alpha)-1}{\alpha}} |\Im(z)|^4}.
\end{align}
Thus, one can write:
\be
G_{N}(z)_{N+k,N+k} = \left[ z -   \sum_{t=1}^{N} r_k(t)^2  \E\left[\hat{G}^{(N+k)}_{N}(z)_{tt}\right] + \epsilon^1_{N,k}(z) + \epsilon^2_{N,k}(z) \right]^{-1}
\ee
where $\epsilon^2_{N,k}(z)$ is a complex valued random variable such that for all $N \in \N$ and $1\leq k \leq N+1$, 
\be
\E[|\epsilon^2_{N,k}(z)|^2] < \frac{C (\ln N)^2}{N^{\frac{\zeta(2 \alpha)-1}{\alpha}} |\Im(z)|^4}.
\ee
In addition, using Lemma \ref{termDown}, we can show:
\begin{align}
&\E\left[\left|  \sum_{t=1}^{N} r_k(t)^2 \left(\E\left[\hat{G}^{(N+k)}_{N}(z)_{tt} - G_{N}(z)_{tt}\right] \right)\right|\right] \\
&\leq  \sum_{t=1}^{N}\E[r_k(t)^2]   \E\left[\left|\hat{G}^{(N+k)}_{N}(z)_{tt} - G_{N}(z)_{tt}\right|\right] \\
&\leq \frac{1}{|\Im(z)|^2\sqrt{N}}.
\end{align}
It follows:
\be
G_{N}(z)_{N+k,N+k} = \left[ z -   \sum_{t=1}^{N} r_k(t)^2  \E\left[G_{N}(z)_{tt}\right] + \epsilon^1_{N,k}(z) + \epsilon^2_{N,k}(z) + \epsilon^3_{N,k}(z) \right]^{-1}
\ee
where $\epsilon^3_{N,k}(z)$ is a complex valued random variable such that for all $N \in \N$ and $1\leq k \leq N+1$, 
\be
\E\left[|\epsilon^3_{N,k}(z)|\right] < \frac{1}{|\Im(z)|^2\sqrt{N}}.
\ee
Let us denote by $I_{N}^t$ the interval $[\frac{t-1}{N},\frac{t}{N}]$. Then we have:
\begin{lemma}
The following inequality holds:
\begin{align*}
\E\left[\left|  \sum_{t=1}^{N}\Big(r_k(t)^2 -M^k(I_{N}^t)\Big)\E\left[G_{N}(z)_{tt}\right]  \right|^2 \right] & \leq \frac{C}{N^{1-\gamma^2}|\Im(z)|^2}
\end{align*}
for some positive constant $C$.
\end{lemma}

\proof We expand the square and, because $r_k(t)$ and $r_{k}(t')$ are independent for $t\not =t'$ conditionally to $M^k$, we have:
\begin{align*}
\E\Big[\Big|  \sum_{t=1}^{N}&\Big(r_k(t)^2 -M^k(I_{N}^t)\Big)\E\left[G_{N}(z)_{tt}\right] \Big|^2 \Big]\\&=\sum_{t,t'=1}^{N}\E\Big[ \Big(r_k(t)^2 -M^k(I_{N}^t)\Big)\Big(r_k(t')^2 -M^k(I_{N}^{t'})\Big)\E\left[G_{N}(z)_{tt}\right] \E\left[G_{N}(z)_{t't'}\right]   \Big]\\
&=\sum_{t=1}^{N}\E\Big[\Big(r_k(t)^2 -M^k(I_{N}^t)\Big)^2\Big]\E\left[G_{N}(z)_{tt}\right] ^2  \\
&=2\sum_{t=1}^{N}\E\left[\left( M^k(I_{N}^t)\right)^2\right]\E\left[G_{N}(z)_{tt}\right] ^2\\
 &\leq 2C\frac{N}{N^{\zeta(2)}|\Im(z)|^2}
\end{align*}\qed

We can thus write
\begin{align}
G_{N}(z)_{N+k,N+k} &= \Bigg[ z -   \sum_{t=1}^{N}  M^k(I_{N}^t) \E\left[G_{N}(z)_{tt}\right] \\
&+ \epsilon^1_{N,k}(z) + \epsilon^2_{N,k}(z) + \epsilon^3_{N,k}(z) + \epsilon^4_{N,k}(z) \Bigg]^{-1}
\end{align}
where $\epsilon^4_{N,k}(z)$ is a complex valued random variable  such that for all $N \in \N$ and $1\leq k \leq N+1$, 
\be
\E\left[|\epsilon^4_{N,k}(z)|^2\right] \leq \frac{C}{N^{\zeta(2)-1}|\Im(z)|^2}.
\ee
Set $\epsilon_{N,k}(z) = \epsilon^1_{N,k}(z) + \epsilon^2_{N,k}(z) + \epsilon^3_{N,k}(z) + \epsilon^4_{N,k}(z)$ and rewrite:
\be
G_{N}(z)_{N+k,N+k} = \left[ z -  \sum_{t=1}^{N} M^k(I_{N}^t)  \E\left[G_{N}(z)_{tt}\right] + \epsilon_{N,k}(z) \right]^{-1}
\ee

We now need to introduce the truncated Radon measure $M^k_\epsilon(dx)$ with Lebesgue density $e^{\omega^k_\epsilon(x)}$ which converges almost surely as $\epsilon$ goes to $0$, in the sense of weak convergence in the space of 
Radon measure, to the measure $M^k$ (see section \ref{cons:loglevy}). 
 
\begin{lemma}\label{epsCV}
For $\epsilon>0$, the following uniform bound holds:
$$ \sup_{N }\E\big[|\sum_{t=1}^{N}  M^k (I_{N}^t)  \E\left[G_{N}(z)_{tt}\right] -\sum_{t=1}^{N}  M^k_\epsilon(I_{N}^t)  \E\left[G_{N}(z)_{tt}\right] |^2\big]\leq \frac{C\epsilon^{1-\gamma^2}}{|\Im(z)|^2}.$$
\end{lemma}

\proof We expand the square. Note that the covariance function $\rho_\epsilon$ of the process $\omega_\epsilon$  increases as $\epsilon$ decreases to $0$ and uniformly  converges as $\epsilon\to 0$ towards 
$\ln_+\frac{\tau}{|x|}$ over the complement of any ball centered at $0$. Thus we have:
\begin{align*}
\sup_{N}\E\big[|\sum_{t=1}^{N} & M^k (I_{N}^t)  \E\left[G_{N}(z)_{tt}\right] -\sum_{t=1}^{N}  M^k_\epsilon(I_{N}^t)  \E\left[G_{N}(z)_{tt}\right] |^2\big]\\
& =\sup_{N}\sum_{t,t'=1}^{N}\E\big[  (M^k (I_{N}^t)  -M^k_\epsilon(I_{N}^t)) (M^k (I_{N}^{t'})  -M^k_\epsilon(I_{N}^{t'}))  \big]\E\left[G_{N}(z)_{tt}\right] \E\left[G_{N}(z)_{t't'}\right] \\
&= \sup_{N}\sum_{t,t'=1}^{N}\E\left[  (M^k (I_{N}^t)  -M^k_\epsilon(I_{N}^t)) (M^k (I_{N}^{t'})  -M^k_\epsilon(I_{N}^{t'}))\right] \E\left[G_{N}(z)_{tt}\right] \E\left[G_{N}(z)_{t't'}\right]   \\
&= \sup_{N}\sum_{t,t'=1}^{N}\left(\E\left[  M^k (I_{N}^t)  M^k (I_{N}^{t'}) \right] -\E\left[  M^k_\epsilon(I_{N}^t)  M^k_\epsilon (I_{N}^{t'}) \right]  \right) \E\left[G_{N}(z)_{tt}\right] \E\left[G_{N}(z)_{t't'}\right] \\
&= \sup_{N}\sum_{t,t'=1}^{N}\E\left[G_{N}(z)_{tt}\right] \E\left[G_{N}(z)_{t't'}\right] \int_{I_{N}^t}\int_{I_{N}^{t'}}\big(e^{\psi(2)\ln_+\frac{\tau}{|r-u|}}-e^{\psi(2)\rho_\epsilon(r-u)}\big)\,drdu\\
&\leq\frac{1}{|\Im(z)|^2}\int_0^{1}\int_0^{1} \big(e^{\psi(2)\ln_+\frac{\tau}{|r-u|}}-e^{\psi(2)\rho_\epsilon(r-u)}\big)\,drdu.
\end{align*}
where, in the fourth line, we used the fact that, if $\mathcal{F}_\epsilon$is the sigma field generated by the random variables
$\mu(A), A\in \mathcal{B}(\{(t,y):y\geq\epsilon\})$, then $\E[M^k(A)|\mathcal{F}_\epsilon] = M_\epsilon^k(A)$ for all borelian set $A$. 
A straightforward computation leads to the relation
\be\label{roeps}\rho_\epsilon(t)=\left\{\begin{array}{ll}\ln\frac{\tau}{\epsilon}+1-\frac{|t|}{\epsilon}&\text{ if }|t|\leq \epsilon\\
\ln\frac{\tau}{|t|}&\text{ if }\epsilon\leq |t|\leq \tau\\
0&\text{ if }\tau<|t|
\end{array}\right.
\ee
By using the expression of $\rho_\epsilon$, it is then plain to obtain the desired bound.\qed

We can thus write
\begin{align}
G_{N}(z)_{N+k,N+k} &= \Bigg[ z -   \sum_{t=1}^{N}  M^k_\epsilon(I_{N}^t) \E\left[G_{N}(z)_{tt}\right] + \epsilon _{N,k}(z) + \delta(\epsilon,N,z) \Bigg]^{-1},
\end{align}
where 
\be\label{limye}
\sup_N\E[|\delta(\epsilon,N,z) |^2]\to 0\quad \text{as }\epsilon\to 0,
\ee
and also:
\be\label{fr}
\E\left[G_{N}(z)_{N+k,N+k}\right] = \E\left[\left[ z -    \sum_{t=1}^{N}  M^k_\epsilon(I_{N}^t)  \E\left[G_{N}(z)_{tt}\right] + \epsilon_{N,k}(z)+ \delta(\epsilon,N,z) \right]^{-1}\right].
\ee
The next step is to study the convergence of the above quantity. Hence we prove (see the proof in the appendix):


\begin{lemma}\label{cvlevy}
The random variable $\sum_{t=1}^{N}  M^k_\epsilon(I_{N}^t)  \E\left[G_{N}(z)_{tt}\right]$ converges in probability as $N\to +\infty$ towards $\int_0^1K_z(x) M^k_\epsilon(dx)$.
\end{lemma}

We fix $\epsilon>0$. For that $\epsilon$, the family of random variables $(\delta(\epsilon,N,z))_N$ is bounded in $L^2$ so that it is tight. Even if it means extracting again a subsequence we assume 
that the couple $(\sum_{t=1}^{N}  M^k_\epsilon(I_{N}^t)  \E\left[G_{N}(z)_{tt}\right],\delta(\epsilon,N,z))_N$ converges in law towards the couple $(\int_0^1K_z(x)M^k_\epsilon(dx),Y_\epsilon)$.  
We remind the reader of \eqref{defgnk} which implies that 
$$\left|\left( z -    \sum_{t=1}^{N}  M^k_\epsilon(I_{N}^t)  \E\left[G_{N}(z)_{tt}\right] + \epsilon_{N,k}(z)+ \delta(\epsilon,N,z) \right)^{-1}\right|\leq \frac{1}{|\Im(z)|}.$$
The   quantity $\left( z -    \sum_{t=1}^{N}  M^k_\epsilon(I_{N}^t)  \E\left[G_{N}(z)_{tt}\right] + \epsilon_{N,k}(z)+ \delta(\epsilon,N,z) \right)^{-1}$ is therefore bounded uniformly with respect to $N,\epsilon$ 
and converges in law towards
$$\left( z -    \int_0^1K_z(x)M^k_\epsilon(dx)+Y_\epsilon\right)^{-1}.$$
We deduce that the expectation of the former quantity converges as $\epsilon\to 0$ towards the expectation of the latter quantity. From \eqref{fr}, we deduce that 
\be
\mu^2_z= \E\left[\left( z - \int_0^1  K_z(x) M^k_\epsilon(dx)+Y_\epsilon \right)^{-1}\right].
\ee
Clearly, standard arguments prove that $\int_0^1  K_z(x) M^k_\epsilon(dx)$ converges almost surely towards $\int_0^1  K_z(x) M^k (dx)$ as $\epsilon\to 0$ ($K_z$ is deterministic (see lemma \ref{concentration}), measurable and bounded) 
and, because of  \eqref{limye}, $Y_\epsilon$ converges almost surely towards $0$ as $\epsilon\to 0$.
Again, because the quantity $\left( z - \int_0^1  K_z(x) M^k_\epsilon(dx)+Y_\epsilon \right)^{-1}$ is bounded uniformly with respect to $\epsilon$, we deduce that:
\be\label{fundeq1}
\mu^2_z= \E\left[\left( z - \int_0^1  K_z(x) M^k (dx)  \right)^{-1}\right].
\ee

\subsection{Second equation}\label{section:secondequation}
Now we turn our attention to the terms $G_{N}(z)_{kk}$ for $k \in \{1,\dots, N\}$.
Again, by using the Schur complement formula, we can write, for $k \in \{1,\dots, N\}$:
\begin{align}
G_{N}(z)_{kk} &= \left[ z - \sum_{i,j=1}^{N} r_i(k) r_j(k) G^{(k)}_{N}(z)_{N+i,N+j} \right]^{-1} \\
                 &= \left[ z - \sum_{i=1}^{N} r_i(k)^2 G^{(k)}_{N}(z)_{N+i,N+i} + \eta^1_{N,k}(z) \right]^{-1}
\end{align}
where, using Lemma \ref{lem:diag2}, $\eta^1_{N,k}(z)$ is a complex valued random variable for which there exists $c>0$ such that for all $N \in \N$ and $1\leq k \leq N, \E[|\eta^1_{N,k}(z)|^2] < c/N$.

With a further use of the Schur complement formula for the term $G^{(k)}_{N}(z)_{N+i,N+i}$, we obtain:
\be
G_{N}(z)_{kk} = \left[ z - \sum_{i=1}^{N} r_i(k)^2 \left[ z - \sum_{s,t \not =k}^{N} r_i(s) r_i(t) G^{(k,N+i)}_{N}(z)_{st} \right]^{-1} + \eta^1_{N,k}(z) \right]^{-1}
\ee
where $G^{(k,N+i)}_{N}(z) = A^{(k,N+i)}_{N}(z)^{-1}$. Note that $G^{(k,N+i)}_{N}(z)$ is independent of $(r_i(t))_{t=1,\dots,N}$. Using the same arguments as in the derivation of the first equation 
(in particular Lemmas \ref{lem:diag1}, \ref{Gchap}, \ref{concsup},  \ref{sup}, \ref{termUP} and  \ref{termDown}), one can show that: 
\begin{align}\label{gnkk}
G_{N}(z)_{kk} = \left[ z -   \sum_{i=1}^{N} \frac{r_i(k)^2}{ z - \sum_{t =1}^{N} M^i(I^t_N) \E\left[G_{N}(z)_{tt}\right] + \delta_{N,k,i}(z)} 
+ \eta^1_{N,k}(z) \right]^{-1}
\end{align}
where $(\delta_{N,k,i}(z))_{1\leq i\leq N}$ are complex random variable such that
\begin{equation}
\E[|\delta_{N,k,i}(z)|]\leq \frac{C}{N^{\min(\frac{1-\gamma^2}{4},\frac{\zeta(2\alpha)-1}{\alpha})}}
\end{equation}
for some positive constant $C$ that does not depend on $i,N$ and for $\alpha > 1$ such that $\zeta(2\alpha)>1$.
\begin{lemma}\label{delta_petit}
One can write:
\begin{align}\label{gnkk2}
G_{N}(z)_{kk} = \left[ z -   \sum_{i=1}^{N} \frac{r_i(k)^2}{ z - \sum_{t =1}^{N} M^i(I^t_N) \E\left[G_{N}(z)_{tt}\right] } 
+ \eta^1_{N,k}(z)+  \eta^2_{N,k}(z) \right]^{-1}
\end{align}
where $ \eta^2_{N,k}(z)$ is a random variable that tends to $0$ in probability as $N$ goes to $\infty$.
\end{lemma}

\proof By using Lemma \ref{lem:im}, we deduce that:
\begin{align}
\sum_{i=1}^{N} \Big|&\frac{r_i(k)^2}{ z - \sum_{t =1}^{N} M^i(I^t_N) \E\left[G_{N}(z)_{tt}\right] + \delta_{N,k,i}(z)} -\frac{r_i(k)^2}{ z - \sum_{t =1}^{N} M^i(I^t_N) \E\left[G_{N}(z)_{tt}\right] } \Big|\nonumber\\
\label{delta}&\leq \frac{1}{|\Im(z)|^2}\sum_{i=1}^{N}r_i(k)^2\min(|\delta_{N,k,i}(z)|,2).
\end{align}
We stress that the lemma is proved as soon as we can prove that   the left-hand side in \eqref{delta} converges in probability to $0$. Hence it is enough 
to prove that $$\E\left[\sum_{i=1}^{N}r_i(k)^2\min(|\delta_{N,k,i}(z)|,2)\right] $$ converges to $0$ as $N$ tends to $\infty$. By noticing 
that:
\begin{equation}
\delta_{N,k,i}(z)= \sum_{s,t \not =k}^{N} r_i(s) r_i(t) G^{(k,N+i)}_{N}(z)_{st} - \sum_{t =1}^{N} M^i(I^t_N) \E\left[G_{N}(z)_{tt}\right], 
\end{equation}
it is straightforward to see that the variables $\Big(r_i(k)^2\min(|\delta_{N,k,i}(z)|,2)\Big)_{1\leq i\leq N+1}$ are identically distributed. Thus we have
$$\E\big[\sum_{i=1}^{N}r_i(k)^2\min(|\delta_{N,k,i}(z)|,2)\big]= N \E\big[ r_1(k)^2\min(|\delta_{N,k,1}(z)|,2)\big].$$ 
Then for all $A>1$ and 
$\alpha>0$, we have
\begin{align*}
N \E\big[ r_1(k)^2\min(|\delta_{N,k,1}(z)|,2)\big] =& N \E\big[ r_1(k)^2\min(|\delta_{N,k,1}(z)|,2)\ind_{\{N r_1(k)^2\leq A\}}\big]\\
&+ N \E\big[ r_1(k)^2\min(|\delta_{N,k,1}(z)|,2)\ind_{\{ N r_1(k)^2> A\}}\big]\\
\leq & A\E[|\delta_{N,k,1}(z)|]+2\E\big[ N r_1(k)^2\ind_{\{N r_1(k)^2> A\}}\big]\\
\leq &\frac{AC}{N^{\frac{\zeta(2)-1}{4}}}+\frac{2}{  A^\alpha}\E\big[N^{1+\alpha}r_1(k)^{2(\alpha+1)} \big]\\
=&\frac{AC}{N^{\frac{\zeta(2)-1}{4}}}+\frac{2 N^{1+\alpha}}{  A^\alpha}\E\big[M^1(0,\frac{1}{N})^{\alpha+1} \big]
\end{align*}
By using the scale invariance property of the measure $M^1$, we have:
\begin{align*}
\E\left[M^1(0,1/N)^{\alpha+1} \right]= &\frac{1}{N^{\zeta(1+\alpha)}}\E\left[ M^1(0,1)^{\alpha+1}\right],
\end{align*}
in such a way that
\begin{equation}
N\E\left[ r_1(k)^2\min(|\delta_{N,k,1}(z)|,2)\right] \leq \frac{AC}{N^{\frac{\zeta(2)-1}{4}}}+2\E\left[ M^1(0,1)^{\alpha+1}\right]\frac{N^{\psi(1+\alpha)}}{A^\alpha}.
\end{equation}
Since $\zeta(2)>5-4\zeta'(1)$ (this inequality is clear with $\zeta(q) = (1+\gamma^2/2)q + q^2 \gamma^2/2$ and is due to our hypotheses of Assumption \ref{condzeta} in the more general case), we can choose $p>0$ such that
\begin{equation}
\frac{\zeta(2)-1}{4}>p>1-\zeta'(1)=\psi'(1).
\end{equation}
The mapping $\alpha\in ]0,+\infty[\mapsto p\alpha-\psi(1+\alpha)$ reduces to $0$ for $\alpha=0$ and, because $p> \psi'(1)$, is strictly positive for $\alpha>0$ small enough. So we choose $\alpha<1$ such 
that $p\alpha-\psi(1+\alpha)>0$ and we set $A=N^p$. We obtain:
$$N \E\big[ r_1(k)^2\min(|\delta_{N,k,1}(z)|,2)\big] \leq \frac{C}{N^{\frac{\zeta(2)-1}{4}-p}}+2^{2+\alpha}\E\big[ M^1(0,T)^{\alpha+1}\big]  \frac{1}{N^{\alpha p-\psi(1+\alpha)}}  .$$
The result follows by letting $N\to \infty$ since $\min((\zeta(2)-1)/4-p,\alpha p-\psi(1+\alpha))>0$.\qed

\begin{lemma}\label{esper}
There exists a constant $c>0$, which does not depend on $N$, such that for each $N \in \N$:
\begin{align*}
\E\Bigg[\bigg| \sum_{i=1}^{N} &\bigg( \frac{r_i(k)^2}{ z -   \sum_{t =1}^{N} M^i(I^t_N) \E\left[G_{N}(z)_{tt}\right]}
- \E\bigg[ \frac{r_i(k)^2}{ z -  \sum_{t =1}^{N} M^i(I^t_N) \E\left[G_{N}(z)_{tt}\right]} \bigg]
 \bigg)  \bigg|^2 \Bigg] \leq \frac{c}{N^{1-\gamma^2}}.
\end{align*}
\end{lemma}
 
\vspace{2mm}
\proof
The proof is straightforward using the fact that for $i\in\{1, \dots,N\}$, the random variables 
\be
\frac{r_i(k)^2}{ z -\sum_{t =1}^{N} M^i(I^t_N) \E\left[G_{N}(z)_{tt}\right]}
\ee 
are i.i.d. random variables and Lemma \ref{lem:im}.
\qed

Therefore we can write
\begin{align}\label{gnkk3}
G_{N}(z)_{kk} = \left[ z -   \sum_{i=1}^{N} \frac{r_i(k)^2}{ z - \sum_{t =1}^{N} M^i(I^t_N) \E\left[G_{N}(z)_{tt}\right] } 
+ \eta^1_{N,k}(z)+  \eta^2_{N,k}(z) +\eta^3_{N,k}(z)\right]^{-1}
\end{align}
with $\E[(\eta^3_{N,k}(z))^2]\leq \frac{c}{N^{1-\gamma^2}}$.

Now we can take the expectation in \eqref{gnkk3} to obtain
\begin{align*}
\E&[L^{1,z}_N(f)]
\\&=\frac{1}{N}\sum_{k=1}^Nf(k/N)\E[G_{N}(z)_{kk}] \\
&=\frac{1}{N}\sum_{k=1}^Nf(k/N)\E \Big[\Big( z -  \E\Big[ \sum_{i=1}^{N} \frac{r_i(k)^2}{ z - \sum_{t =1}^{N} M^i(I^t_N) \E \big[G_{N}(z)_{tt} \big] } 
\Big]+ \eta_{N,k}(z)\Big)^{-1}\Big] \\
&= \frac{1}{N}\sum_{k=1}^Nf(k/N)\E \Big[\Big( z -  N \E\Big[ \frac{M\left[\frac{k-1}{N};\frac{k}{N}\right]}{ z - \sum_{t =1}^{N} M(I^t_N) \E \big[G_{N}(z)_{tt} \big] } 
\Big]+ \eta_{N,k}(z) \Big)^{-1}\Big]
\end{align*}
with $\eta_{N,k}(z) = \eta^1_{N,k}(z)+\eta^2_{N}(z)+\eta^3_{N,k}(z)$.
Then, by introducing the truncated measure $M_\epsilon$ and by using the Girsanov formula, we can approximate (uniformly in $N$) this last expression by:
\begin{equation}
\frac{1}{N}\sum_{k=1}^Nf(k/N)\E \Big[\Big( z -  N \E\Big[ \frac{M_\epsilon\left[\frac{k-1}{N};\frac{k}{N}\right]}{ z - \sum_{t =1}^{N} M_\epsilon(I^t_N) \E \big[G_{N}(z)_{tt} \big] } 
\Big] \Big)^{-1} + \hat{\delta}(N,k,z,\epsilon) \Big]
\end{equation}
with $\sup_{N,k} \E[|\hat{\delta}(N,k,z,\epsilon)|^2]$ going to $0$ when $\epsilon$ is going to $0$.
Along some appropriate subsequence, this latter quantity converges as $N \to + \infty$ to:
\begin{equation}
\int_0^1f(x)\E\left[\left( z -   \E\left[  \frac{e^{\omega_\epsilon(x)}}{ z -\int_0^1   K_z(r)\,M_\epsilon(dr) }\right] \right)^{-1} + Y^\epsilon\right]\,dx
\end{equation}
where $Y^\epsilon$ is such that $\E[(Y^\epsilon)^2]$ converges to $0$ when $\epsilon$ is going to $0$. 
And, we thus obtain gathering the above arguments that:
\be\label{eps2}
\begin{split}
\int_0^1f(x)K_z(x)\,dx = \int_0^1f(x)\E\left[\left( z -   \E\left[  \frac{e^{\omega_\epsilon(x)}}{ z -\int_0^1   K_z(r)\,M_\epsilon(dr) }\right] \right)^{-1} + Y^\epsilon\right]\,dx.
\end{split}
\ee
It remains to pass to the limit as $\epsilon\to 0$ in that expression. This job is carried out with the help of a Girsanov type transform in Appendix \ref{girs}.\qed

\subsection{Uniqueness of the solution to the system of equations}\label{section:unicity}

Let $X$ be the space of bounded measurable functions $[0,1] \to \C$ endowed with the uniform norm defined for $f \in X$ by:
\begin{equation}
||f||_\infty = \sup_{x \in [0,1]} |f(x)|.
\end{equation}
Define the operator $T:X \to X$ by setting, for $g \in X$ and for all $x\in[0,1]$:
\begin{equation}
Tg(x) = \frac{1}{z- q\E\left[\left( z - \int_0^1 \left(\frac{\tau}{|t-x|}\right)_+^{\gamma^2} g(t) M(dt) \right)^{-1}\right]} 
\end{equation}
For $g,h \in X$ and for all $x\in[0,1]$, we have:
\begin{align*}
|Tg(x) - Th(x)| &\leq \frac{q}{|\Im(z)|^4} \E\left[\int_0^1 \left(\frac{\tau}{|t-x|}\right)_+^{\gamma^2} |g(t) - h(t)|  M(dt) \right] \\
                &\leq \frac{q}{|\Im(z)|^4} \E\left[\int_0^1 \left(\frac{\tau}{|t-x|}\right)_+^{\gamma^2}  M(dt) \right] ||g-h||_{\infty} \\
                &\leq \frac{q}{|\Im(z)|^4} \int_0^1 \left(\frac{\tau}{|t-x|}\right)_+^{\gamma^2} dt ||g-h||_{\infty}. 
\end{align*}
 
Recall that $\gamma^2 <1/3$, and thus it is easy to see that:
\begin{equation}
\sup_{x \in [0,1]} \int_0^1 \left(\frac{\tau}{|t-x|}\right)_+^{\gamma^2} dt < +\infty
\end{equation} 
And we can deduce that there exists a positive constant $C$ such that:
\begin{align}
\sup_{x \in [0,1]} |Tg(x) - Th(x)| \leq \frac{C}{|\Im(z)|^4} ||g-h||_{\infty} 
\end{align}
If $z$ is such that $C/|\Im(z)|^4 < 1$, the operator $T$ is contracting and thus has a unique fixed point $g$ in the Banach space $X$. 
We conclude that, for each $z$ with $|\Im(z)|$ large enough, there exists a unique bounded function $K_z :[0,1] \to \C$ such that for all $x\in[0,1]$:
\begin{equation}
K_z(x) = \frac{1}{z- q\E\left[\left( z - \int_0^1 \left(\frac{\tau}{|t-x|}\right)_+^{\gamma^2} K_z(t) M(dt) \right)^{-1}\right]}. 
\end{equation}

Using the first equation, it is now plain to see that, for $z$ such that $C/|\Im(z)|^4 < 1$, the constant $\mu_z^2$ is uniquely defined by the system of equations 
(by the first equation, it is a function of the function $K_z$, which is uniquely defined for such $z$).  

Now it remains to show that the limit point $\mu_z^2$ is uniquely defined for all $z\in \C \setminus \R$. It will be easy to see using analyticity arguments.
Indeed, from the Montel theorem, every limit point $\mu^2_z$ is holomorphic on the set $\C \setminus \R$ since it is the pointwise limit of a subsequence of the sequence of holomorphic functions 
$L^{1,z}_N([0,1])$ that are uniformly bounded on each compact set of $\C\setminus\R$ (see Lemma \ref{lem:im}). Thus, $\mu_z^2$ is uniquely defined for each $z\in \C \setminus \R$ by analytic extension
(we have just seen that $\mu_z^2$ is uniquely defined for a set of $z$ with accumulation points).  

The same argument holds for the unicity of the integral $\int_0^1 K_z(x) dx$. Indeed, every limit point $\int_0^1 K_z(x)dx$ is a holomorphic function on $\C \setminus \R$ that has some prescribed value on the set
$\{z \in \C \setminus \R: C/|\Im(z)|^4 < 1 \}$, which has accumulation points.

\subsection{Proof of Theorem \ref{main2}, \ref{main1} and \ref{main3}}
Let us gather the above arguments to prove the main theorems. 

\emph{Proof of theorem \ref{main3}}: 
it is a direct consequence of sections \ref{section:firstequation}, \ref{section:secondequation} and \ref{section:unicity}.

\emph{Proof of theorem \ref{main1} i)}: 
The limit points $K_z(x) dx$ and $\mu_z^2 dx$ of the two complex measures $\E[L^{1,z}_N]$ and $\E[L^{2,z}_N]$ are uniquely defined because 
$\mu_z^2$ and $K_z(x)$ satisfy a fixed point system of equations (we have just seen this in theorem  \ref{main3}).

\emph{Proof of theorem \ref{main1} iii)}:
We need to prove that $\mu^2_z$ is the Stieltjes transform of a probability measure $\upsilon$. From \cite{Geronimo}, it suffices to prove that $\mu^2_z$ is holomorphic 
over $\C\setminus\R$, maps $\{z\in\C\setminus\R; \Im(z)<0\}$ to $\{z\in\C\setminus\R; \Im(z)>0\}$ and that $\lim_{y\to \infty} iy \mu^2_{iy} = 1$ ($y\in\R$). Let us check those properties.  
We have already seen in section \ref{section:unicity} that $\mu_z^2$ is holomorphic. From Lemma \ref{lem:im}, $\mu^2$ maps $\{z\in\C\setminus\R; \Im(z)<0\}$ to $\{z\in\C\setminus\R; \Im(z)>0\}$.
Finally, from Theorem \ref{main3}, we have
$$z\mu^2_z=\E\left[\frac{1}{1-z^{-1}\int_0^1K_z(x)\,M(dx)}\right].$$ As $|K_z(x)|\leq |\Im(z)|^{-1}$, the term $\int_0^1K_z(x)\,M(dx)/z$ converges pointwise towards $0$ when $z=iy$ and $y\to \infty$. 
Furthermore, from Lemma \ref{imK}, we have $\Im(z)\Im(K_z(x)\leq 0$ in such a way that $\left|z-\int_0^1K_z(x)\,M(dx)\right|^{-1} \leq  |\Im(z)|^{-1}$. Therefore
$$\big|\frac{z}{z-\int_0^1K_z(x)\,M(dx)}\big|\leq 1$$ when $z$ takes on the form $z=iy$ ($y\in\R$). The dominated convergence theorem then implies that $\lim_{y\to \infty}iy\mu^2_{iy}=1$ and we can conclude $\mu^2$ 
is indeed the Stieltjes 
transform of a (unique) probability measure $\upsilon$.

\emph{Proof of theorem \ref{main2} i) and \ref{main1} ii)}
We observe that, for $z \in \C \setminus \R$:
\begin{equation}\label{eq:decomp}
A_N(z) 
\begin{pmatrix}
         z I_{T} & 0 \\
         X_N & z I_{N}
\end{pmatrix} 
=
\begin{pmatrix}
         z^2 I_{T} -   {}^{t} \! X_N X_N & - z  {}^{t} \!  X_N  \\
         0 & z^2 I_{N}
\end{pmatrix} .
\end{equation}

Let us rewrite the matrix $G_N(z)=A_N(z)^{-1}$ under the form: 
\begin{equation}
G_N(z)=\begin{pmatrix}
         G_1(z) & {}^{t} \! G_{1,2}(z) \\
         G_{1,2}(z) & G_2(z)
\end{pmatrix},
\end{equation}
where $G_1(z),G_{1,2}(z),G_2(z) $ are respectively of size $T\times T$, $N\times T$, $N\times N$.

By taking the inverse in the relation \eqref{eq:decomp}, we obtain:
\begin{equation}\label{dectamer}
\begin{pmatrix}
          I_{T}/z & 0 \\
         -X_N / z^2 &  I_{N}/z
\end{pmatrix} 
\begin{pmatrix}
         G_1(z) & {}^{t} \! G_{1,2}(z) \\
         G_{1,2}(z) & G_2(z)
\end{pmatrix} 
=
\begin{pmatrix}
         (z^2 I_{T} -   {}^{t} \! X_N X_N)^{-1} & B \\
         0 & I_{N} / z^2
\end{pmatrix} 
\end{equation}
where $B=(z^2 I_{T} -  {}^{t} \! X_N X_N)^{-1} {}^{t} \!  X_N / z$.

It can be rewritten, using the fact that $-X_N G_1(z)+ z G_{1,2}(z) = 0$ and $-X_N {}^{t} \! G_{1,2}(z) + z G_2(z) = I_N$, as:
\begin{equation}\label{dectamer}
\begin{pmatrix}
         G_1(z)/z &  {}^{t} \! G_{1,2}(z)/z \\
         0 & I_N/z^2
\end{pmatrix} 
=
\begin{pmatrix}
         (z^2 I_{T} -   {}^{t} \! X_N X_N)^{-1} & B \\
         0 & I_{N} / z^2
\end{pmatrix} 
\end{equation}

Therefore, taking the trace we get: 
\begin{equation}\label{eq:trace}
\frac{1}{T z}\sum_{k=1}^{T} G_N(z)_{kk}= \frac{1}{T} \text{tr} (z^2 I_{T} -  {}^{t} \! X_N X_N)^{-1},
\end{equation}

and, by using the fact that the eigenvalues of ${}^{t} \! X_N X_N$ are those of $X_N  {}^{t} \! X_N$ augmented with $T-N$ zeros:
\begin{equation}
\frac{1}{T z} \sum_{k=1}^{T} G_N(z)_{kk} = \frac{1}{T} \text{tr} (z^2 I_{N} -  X_N {}^{t} \!  X_N)^{-1} + \frac{T-N}{T z^2}. 
\end{equation}

Now, taking expectation and using theorem \ref{main1}, we deduce:
\begin{equation}\label{eqA}
\int_0^1 K_z(x) dx = q z \lim_{N\to \infty} \frac{1}{N} \E\left[\text{tr} (z^2 I_{N} -  X_N {}^{t} \!  X_N)^{-1}\right] + \frac{1-q}{z}
\end{equation}

Using the fact that (by (\ref{rel:rn})) the spectrum of $B_N$ contains $2N$ eigenvalues which are the positive and negative square-roots of the spectrum of $R_N={}^{t} \! X_N X_N$ plus $T-N$ zero 
eigenvalues and the fact that 
$1/(z-\lambda)+1/(z+\lambda)=2z/(z^2-\lambda^2)$, we can see that:
\begin{equation}\label{eq:as1}
\frac{1}{N+T} \sum_{k=1}^{N+T} G_N(z)_{kk}= \frac{2z}{N+T} \text{tr} (z^2 I_{N} -  X_N {}^{t} \! X_N)^{-1}+\frac{T-N}{T+N} \frac{1}{z}
\end{equation}

Using the relation \ref{rel:trace2} and theorem \ref{main1}, it is easy to see that:
\begin{equation}\label{eq:TRdeG_N}
\lim_{N\to +\infty} \frac{1}{N+T} \sum_{k=1}^{N+T} \E[G_N(z)_{kk}] = \frac{1}{1+q} \left(q \mu_z^2 + \int_0^1 K_z(x) dx \right) 
\end{equation}

Taking expectation in \ref{eq:as1} and using (\ref{eq:TRdeG_N}), we get:
\begin{align}\label{eqB}
\frac{1}{1+q} \left(q \mu_z^2 + \int_0^1 K_z(x) dx \right) &= \frac{2q z}{1+q} \lim_{N\to \infty} \frac{1}{N} \E\left[\text{tr} (z^2 I_{N} -  X_N {}^{t} \!  X_N)^{-1}\right] \\ 
&+ \frac{1-q}{1+q} \frac{1}{z}.
\end{align}

From equations (\ref{eqA}) and (\ref{eqB}), we get the following relation:
\begin{equation}\label{equationKmu}
\int_0^1 K_z(x) dx = q \mu_z^2 + \frac{1-q}{z}. 
\end{equation}
and theorem \ref{main1} ii). is proved. 

With (\ref{equationKmu}), (\ref{eq:TRdeG_N}) becomes:
\begin{equation}\label{eq:TRdeG_N2}
\lim_{N\to +\infty} \frac{1}{N+T} \sum_{k=1}^{N+T} \E[G_N(z)_{kk}] = \frac{1}{1+q} \left(2 q \mu_z^2 + \frac{1-q}{z} \right) 
\end{equation} 
and, we note that the right hand side of (\ref{eq:TRdeG_N2}) is the Stieltjes transform of the measure $2q/(1+q) \upsilon(dx) + (1-q)/(1+q) \delta_0(dx)$. 
Thus, the mean spectral measure $\E[\mu_{B_N}]$ converges weakly to the measure  $2q/(1+q) \upsilon(dx) + (1-q)/(1+q) \delta_0(dx)$. 

We have also:
\begin{equation}
\lim_{N\to \infty} \frac{1}{N} \E\left[\text{tr} (z^2 I_{N} -  X_N {}^{t} \!  X_N)^{-1}\right] = \frac{\mu_z^2}{z}
\end{equation}

Again using the fact that, for all $x\in \R, 1/(z^2-x^2)=(1/(z-x)+1/(z+x))/(2z)$ and the fact that $\upsilon(dx)$ is a symmetric measure on $\R$ ($\upsilon(dx)$ is the weak limit 
of $\E\left[\mu_{B_N}\right]$, which is symmetric since the spectrum of $B_N$ is symmetric with respect to $0$ almost surely), we see that:
\begin{align}
\lim_{N\to \infty} \frac{1}{N} \E\left[\text{tr} (z^2 I_{N} -  X_N {}^{t} \!  X_N)^{-1}\right] &= \frac{1}{z} \int_\R \frac{\upsilon(dx)}{z-x} \\
&= \int_\R \frac{\upsilon\circ (x^2)^{-1}(dx)}{z^2-x}. 
\end{align}

This implies that, for each $z \in \C\setminus\R$, 
\begin{align}
\lim_{N\to \infty} \frac{1}{N} \E\left[\text{tr} (z I_{N} -  X_N {}^{t} \!  X_N)^{-1}\right] = \int_\R \frac{\upsilon\circ (x^2)^{-1}(dx)}{z-x}. 
\end{align}
and thus, the probability measure $\E[\mu_{R_N}]$ converges weakly to the measure $\upsilon\circ (x^2)^{-1}(dx)$.

\emph{Proof of theorem \ref{main2} ii)}: using relation (\ref{rel:trace2}) and lemma \ref{concentration}, it is plain to check that $\int_\R (z-x)^{-1} \mu_{B_N}(dx)$
converges in probability to the Stieltjes transform of the probability measure $2q/(1+q)\upsilon(dx)+(1-q)/(1+q)\delta_0(dx)$. This convergence holds for finite dimensional vectors 
$(\int_\R  (z_i - x)^{-1} \mu_{B_N}(dx)), i=1, \dots, d)$ as well. Using the fact that the set of functions $\{(z-x)^{-1}, z \in \C\setminus\R\}$ is dense in the set $C_0(\R)$ of continuous 
functions on $\R$ going to $0$ at infinity, we can show, for each $f \in C_0(\R)$, that $\int f(x) \mu_{B_N}(dx)$ converges in probability to $\int f(x) (2q/(1+q)\upsilon(dx)+(1-q)/(1+q)\delta_0(dx))$. 
But, since $\mu_{B_N}(\R) = 2q/(1+q)\upsilon(\R)+(1-q)/(1+q)\delta_0(\R) = 1$, this vague convergence can be strengthened in a weak convergence. 
With the relations $\mu_{B_N^2} = 2N/(N+T) \mu_{R_N} + (T-N)/(T+N) \delta_0$ and the fact that $\int f(x) \mu_{B_N^2}(dx) = \int f(x^2) \mu_{B_N}(dx)$, it is plain to conclude that
$\mu_{R_N}$ converges weakly in probability to $\upsilon\circ (x^2)^{-1}(dx)$.

\emph{Proof of theorem \ref{main2} iii)}: again using relation (\ref{rel:trace2}) and lemma \ref{concentration} together with Borel-Cantelli's lemma, one can show that the two spectral measures
$\mu_{B_{N_k}}$  converges weakly almost surely to $2q/(1+q)\upsilon(dx)+(1-q)/(1+q)\delta_0(dx)$. It is then easy to deduce as before that $\mu_{R_{N_k}}$ converges weakly almost surely to $\upsilon\circ (x^2)^{-1}(dx)$.


\appendix

\section{Auxiliary lemmas}
\begin{lemma}\label{SommeLigne}
Let $A$ be a $n\times n$ complex matrix such that the Hermitian matrix $M= A\bar{A}^T$ has spectral radius $\lambda_{max}$.
Then, for all $i$, we have:
\begin{equation}
\sum_{j=1}^{n} \mid A_{ij} \mid^2  \leq  \lambda_{max}.
\end{equation}
\end{lemma}

\proof It is straightforward to see that all the entries of $M$ are, in modulus, smaller than $\lambda_{max}$.
On the other hand, we have:
\begin{align*}
M_{ii} = \sum_{j=1}^{n} \mid A_{ij} \mid^2.
\end{align*} 
and, thus:
\begin{equation}
\sum_{j=1}^{n} \mid A_{ij} \mid^2  \leq  \lambda_{max}.
\end{equation}
\hfill$\square$

\begin{lemma}\label{lem:diag1}
There exists $C>0$ such that for each $N \in \N$ and $k \in \{1, \dots, N \}$:
\begin{equation*}
\E\left[ \left| \sum_{s \not = t}^{N}  r_k(s) r_k(t) G^{(N+k)}_{N}(z)_{st}   \right|^2   \right]  \leq  \frac{C}{N^{1-\gamma^2}}.
\end{equation*}
Similarly, for each $N \in \N$ and $k \in \{1, \dots, N \}$, $i \in \{1, \dots, N \}$, we have the following inequality concerning the conditional expectation with respect to $M^i$:
\begin{equation*}
\E\left[ \left| \sum_{s,t \not =k,s\not =t}^{N} r_i(s) r_i(t) G^{(k,N+i)}_{N}(z)_{st}   \right|^2 |M^i  \right]  \leq  \frac{C}{N^{1-\gamma^2}}.
\end{equation*}

\end{lemma}

 \proof We first expand the square and use the independence of $(r_k(s))_s$ from $G^{(N+k)}_{N}(z)$:
\begin{align*}
 \E\left[ \left| \sum_{s \not = t}^{N}  r_k(s) r_k(t) G^{(N+k)}_{N}(z)_{st}   \right|^2   \right] & = 2 \sum_{s \not = t}^{N} \E\left[ r_k(s)^2 r_k(t)^2 \right] \E\left[ \left|  G^{(N+k)}_{N}(z)_{st}   \right|^2   \right] 
 \end{align*}
 Now we compute
  \begin{align*}
 \E\left[ r_k(s)^2 r_k(t)^2 \right]&= \E\left[ M^k(\frac{s-1}{N},\frac{s}{N})M^k(\frac{t-1}{N},\frac{t}{N}) \right]\\
 &=\int_{\frac{s-1}{N}}^{\frac{s}{N}}\int_{\frac{t-1}{N}}^{\frac{t}{N}}\max\left(1,\frac{\tau}{|r-u|}\right)^{\psi(2)}\,drdu\\
 &\leq \int_{0}^{\frac{1}{N}}\int_{\frac{ 1}{N}}^{\frac{2}{N}}\max\left(1,\frac{\tau}{|r-u|}\right)^{\psi(2)}\,drdu
  \end{align*}
  We consider $N$ large enough so as to make $2/N\leq \tau$. The above integral is then plain to compute and we get
\begin{equation} 
\E\left[ r_k(s)^2 r_k(t)^2 \right]\leq \frac{\tau^{\psi(2)}(2^{2-\psi(2)}-2)}{(1-\psi(2))(2-\psi(2))}\frac{1}{N^{2-\psi(2)}}.
\end{equation}
Thus we have for some positive constant $C$
 \begin{align*}
 \E\left[ \left| \sum_{s \not = t}^{N}  r_k(s) r_k(t) G^{(N+k)}_{N}(z)_{st}   \right|^2   \right] 
& \leq \frac{C}{N^{2-\psi(2)}} \sum_{s \not = t}^{N}     \E\left[ \left|  G^{(N+k)}_{N}(z)_{st}    \right|^2   \right]\\
& \leq \frac{C}{N^{1-\psi(2)}}\frac{1}{| \Im(z)|^2}, 
\end{align*}
where we have used the fact that almost surely: 
\begin{equation*}
\frac{1}{2N-1}\sum_{s,t\not = N+k}^{2N}  \left|  G^{(N+k)}_{N}(z)_{st}  \right|^2 \leq \frac{1}{| \Im(z)|^2}.
\end{equation*}
It just remains to see that $\psi(2)=\gamma^2$.  
To prove the second relation, we follow the same argument by noticing that $(r_i(t))_{t}$ and  $G^{(k,N+i)}_{N}(z) $ are independent conditionally to $M^i$.
\qed

 \begin{lemma}\label{lem:diag2}
 There exists some  constant $c>0$ such that for each $N \in \N$ and $k \in \{1, \dots, N \}$:
 \begin{equation*}
 \E\left[ \left| \sum_{i \not = j}^{N}  r_i(k) r_j(k) G^{(k)}_{N}(z)_{N+i,N+j} \right|^2 \right]  \leq  \frac{c}{N}.  
 \end{equation*}
 \end{lemma}

\proof Again we expand the square and we use the fact that, conditionally to the $(M^i)_i$, the quantities $r_i(k),r_j(k), G^{(k)}_{N}(z)_{N+i,N+j}$ are independent and $r_i(k),r_j(k)$ are centered. Indeed, 
conditionally to the $(M^i)_i$, the variables  $r_i(k),r_j(k), G^{(k)}_{N}(z)_{N+i,N+j}$ involve different increments of the Brownian motion. Thus we have
\begin{align*}
\E\left[ \left| \sum_{i \not = j}^{N}  r_i(k) r_j(k) G^{(k)}_{N}(z)_{N+i,N+j} \right|^2 \right]  & =  \sum_{i \not = j}^{N} \E\left[ r_i(k)^2 r_j(k)^2 \right]\E\left[ \left|  G^{(k)}_{N}(z)_{N+i,N+j}  \right|^2 \right] \\
& \leq   \sum_{i \not = j}^{N} \E [ r_i(k)^2] \E[r_j(k)^2] \E\left[  \left|  G^{(k)}_{N}(z)_{N+i,N+j}  \right|^2\right]   \\
& =  N^{-2} \sum_{i \not = j}^{N} \E\left[  \left|  G^{(k)}_{N}(z)_{N+i,N+j}  \right|^2\right]   \\
& \leq  \frac{c}{N},
\end{align*}
where we have used the fact that almost surely: 
\begin{equation*}
\frac{1}{2N-1}\sum_{i,j\not = k}^{2N}  \left|  G^{(k)}_{N}(z)_{i,j}  \right|^2 \leq \frac{1}{| \Im(z)|^2}.
\end{equation*}
\qed

{\it Proof of Lemma \ref{cvlevy}.} We define the function $f^{k,\epsilon}_N$  on the interval $[0,1]$ by 
$$f^{k,\epsilon}_N(x)=NM^{k,\epsilon}(I^t_N) \text{ if } x\in I^t_N.$$ 
Notice the relation:
$$\sum_{t=1}^{N}  M^k_\epsilon(I_{N}^t)  \E\left[G_{N}(z)_{tt}\right]=\int_0^1f^{k,\epsilon}_N(r)\,d\E[L^{1,z}_N](dr).$$
Then, by stationarity, we have:
\begin{align*}
&\E\left[\left|\int_0^1 f^{k,\epsilon}_N(r)\,d\E[L^{1,z}_N](dr)-\int_0^1e^{\omega^k_\epsilon(r)}\,d\E[L^{1,z}_N](dr)\right|\right]\\
&\leq \sum_{t=1}^{N}\E\left[\left|\int_{I^t_N}(f^{k,\epsilon}_N(r)-e^{\omega^k_\epsilon(r)})\,d\E[L^{1,z}_N](dr) \right|\right] \\
&\leq \frac{N}{|\Im(z)|} \sup_{r\in I^1_N} \E\left[\left|\int_{I^1_N}(e^{\omega^k_\epsilon(u)}-e^{\omega^k_\epsilon(r)})\,du \right|\right]\\
&\leq \frac{N}{|\Im(z)|} \sup_{r\in I^1_N} \int_{I^1_N}\E\left[\left|e^{\omega^k_\epsilon(u)}-e^{\omega^k_\epsilon(r)})\right|^2\right]^{1/2}\,du \\
&\leq \frac{N}{|\Im(z)|} \sup_{r\in I^1_N} \int_{I^1_N}\left( 2e^{\psi(2)\rho_\epsilon(0)}-2e^{\psi(2)\rho_\epsilon(r-u)} \right)^{1/2}\,du. 
\end{align*}
Because of the   continuity of the function $\rho_\epsilon$ over $[0,1]$, we have
\be\label{cvproba1}
\E\Big[\Big|\int_0^1f^{k,\epsilon}_N(r)\,d\E[L^{1,z}_N](dr)-\int_0^1e^{\omega^k_\epsilon(r)}\,d\E[L^{1,z}_N](dr)\Big|\Big]\to 0\quad \text{ as }N\to \infty.
\ee
In a quite similar way, we can prove that 
\be\label{cvproba2}
\E\Big[\Big|\int_0^1 e^{\omega^k_\epsilon}*\phi_p(r)\,d\E[L^{1,z}_N](dr)-\int_0^1e^{\omega^k_\epsilon(r)}\,d\E[L^{1,z}_N](dr)\Big|\Big]\to 0\quad \text{ as }p\to \infty \text{ uniformly w.r.t. N}
\ee and 
\be\label{cvproba3}
\E\Big[\Big|\int_0^1 e^{\omega^k_\epsilon}*\phi_p(r)K_z(r)\,dr-\int_0^1e^{\omega^k_\epsilon(r)}K_z(r)\,dr\Big|\Big]\to 0\quad \text{ as }p\to \infty \text{ uniformly w.r.t. N}
\ee
where $(\phi_p)_{p\in\N}$ is a regularizing sequence and $*$ stands for the convolution. Furthermore, for each fixed $p$ and because of the weak convergence of $\E[L^{1,z}_N]$ towards $K_z(x)dx$, we have almost surely
\be\label{cvproba4}
 \int_0^1 e^{\omega^k_\epsilon}*\phi_p(r)\,d\E[L^{1,z}_N](dr)\to \int_0^1 e^{\omega^k_\epsilon}*\phi_p(r)K_z(r)\, dr\quad \text{ as }N\to \infty.
\ee
We prove the result by gathering \eqref{cvproba1}  \eqref{cvproba2}  \eqref{cvproba3} and  \eqref{cvproba4}.\qed

\section{Sup of MRW}
 
 Here we prove 
 \begin{proposition}\label{sup}
 We have for all $k=1,\dots,N+1$
 $$\E\left[\sup_{t=1,\dots,N}r_k(t)^4\right]\leq  C\frac{(\ln N)^2}{N^{\frac{\zeta(2\alpha)-1}{\alpha}}}.$$ for some positive constant $C$.
 \end{proposition}
 
 \proof To prove the result, we first prove 
  \begin{lemma}\label{GS}
There exists a constant $C$ such that,
if $(X_i)_{1\leq i \leq N}$ are iid centered Gaussian random variables then:
\begin{equation*}
\E\left[\max_{1\leq i \leq N}|X_i|^4\right]\leq C \max_{1\leq i \leq N} \E[X_i^2]^2 (\ln N)^2.
\end{equation*}
\end{lemma}
 
\proof By homogeneity, it suffices to assume that $ \E[X_i^2]=1$. Then we have for all $\delta\geq 0$
\begin{align*}
\E\left[\max_{1\leq i \leq N}|X_i|^4 \right] &\leq \delta+N \int_\delta^\infty \P(|X_1|^4>t) dt \\
&\leq  \delta+2N\int_\delta^\infty \P(X_1>t^{1/4}) dt\\
&\leq \delta+\frac{2N}{\sqrt{2 \pi}} \int_\delta^\infty e^{-\sqrt{t}} dt\\
&\leq \delta+\frac{4N}{\sqrt{2 \pi}} \int_{\sqrt{\delta}}^\infty e^{-t} t dt \\
&\leq \delta+\frac{4N}{\sqrt{2 \pi}} \left( \sqrt{\delta} e^{-\sqrt{\delta}} + e^{-\sqrt{\delta}} \right), 
\end{align*}
and this last expression can be made smaller than $C(\ln N)^2$ by choosing $\delta=(\ln N)^2$.\qed
 
 We want apply the above lemma after conditioning with respect to the law of the MRM $M^k$:
 \begin{align*}
 \E\left[\sup_{t=1,\dots,N}r_k(t)^4\right]&= \E\left[\E\left[\sup_{t=1,\dots,N}r_k(t)^4|M^k\right]\right].
 \end{align*}
 Notice then that, conditionally to $M^k(0,\frac{1}{N})=x_1,\dots,M^k(\frac{N-1}{N},1)=x_N$, the vector $(r_k(1),\dots,r_k(N)$ has the same law as the increments 
 of $B$: $(B_{x_1}-B_0,\dots,B_{x_N}-B_{x_{N-1}})$. By applying Lemma \ref{GS}, we deduce that
 $$\E\left[\sup_{t=1,\dots,N}r_k(t)^4|M^k\right]\leq C(\ln N)^2 \max_{t=1,\dots,N}M^k\big(\frac{t-1}{N},\frac{t}{N}\big)^2.$$
Thus we deduce
\begin{equation}\label{step1}
\E\left[\sup_{t=1,\dots,N}r_k(t)^4\right]\leq C(\ln N)^2 \E\left[\left(\max_{t=1,\dots,N}M^k\big(\frac{t-1}{N},\frac{t}{N}\big)\right)^2\right].
\end{equation}
Finally we have for all $\delta>0$ and for $\alpha>1$ such that $\zeta(2\alpha)>1$:
\begin{align*}
\E\left[\left(\max_{t=1,\dots,N}M^k\big(\frac{t-1}{N},\frac{t}{N}\big)\right)^2\right]&\leq \delta+N\int_\delta^\infty\P\big(M^k\big(\frac{t-1}{N},\frac{t}{N}\big)^2>x\big)\,dx\\
&\leq \delta+N\int_\delta^\infty\frac{1}{x^\alpha}\E\big[M^k\big(\frac{t-1}{N},\frac{t}{N}\big)^{2\alpha}\big]\,dx\\
&\leq \delta+C\delta^{1-\alpha}N^{1-\zeta(2\alpha)}
\end{align*}
for some constant $C$ only depending on $\alpha,\tau$ and $\gamma^2$. Choose now $\delta=N^{\frac{1-\zeta(2\alpha)}{\alpha}}$ so as to get
\begin{equation}\label{step2}
\E\left[\sup_{t=1,\dots,N}r_k(t)^4\right]\leq (1+C)\frac{(\ln N)^2}{N^{\frac{\zeta(2\alpha)-1}{\alpha}}}
\end{equation}
\qed
\section{Girsanov transform}\label{girs}

\begin{lemma}\label{lem:girs}
Let $\mu$ be an independently scattered infinitely divisible random measure associated to $(\psi,\theta)$, where 
$$\forall q \in\R,\quad \psi(q)=mq+\frac{1}{2}\sigma^2q^2+\int_\R(e^{qz}-1)\nu(dz), $$
 $\psi(2)<+\infty$ and   $\psi(1)=0$. Let $B$ be a bounded Borelian set. We define a new probability measure $\P_B$ (with expectation $\E_B$) by:
$$\forall A \text{measurable set},\quad \P_B(A)=\E[\ind_Ae^{\mu(B)}] .$$
Then, under $\P_B$, $\mu$ has the same law as $\mu+\mu_B $ where $\mu_B$ is an independently scattered infinitely divisible random measures independent of $\mu$ and is associated to $(\psi_B,\theta_B)$ given by
\begin{align*}
\psi_B(q)&= q\sigma^2+\int_\R(e^{qx}-1)(e^x-1)\nu(dx)\\
\theta_B(\cdot)&=\theta(\cdot\cap B).
\end{align*}
\end{lemma}

\proof It suffices to compute the joint distribution of $p$ disjoint sets $A_1,\dots,A_p$. We have for any $\lambda_1,\dots,\lambda_p\in\R$:
\begin{align*}
\E_B\Big[e^{\lambda_1\mu(A_1)+\dots+\lambda_p\mu(A_p)}\Big]&=\E\Big[e^{\lambda_1\mu(A_1)+\dots+\lambda_p\mu(A_p)+\mu(B)}\Big]\\
&=\E\Big[e^{\lambda_1\mu(A_1\setminus B)+\dots+\lambda_p\mu(A_p\setminus B)+\lambda_1\mu(A_1\cap B)+\dots+\lambda_p\mu(A_p\cap B)+\mu(B)}\Big]\\
&=\E\Big[e^{\lambda_1\mu(A_1\setminus B)+\dots+\lambda_p\mu(A_p\setminus B)+(\lambda_1+1)\mu(A_1\cap B)+\dots+(\lambda_p+1)\mu(A_p\cap B)+\mu(B\setminus\bigcup_{i=1}^nA_i)}\Big]\\
&=\E\Big[e^{\lambda_1\mu(A_1\setminus B)+\dots+\lambda_p\mu(A_p\setminus B)}\Big]\E\Big[e^{(\lambda_1+1)\mu(A_1\cap B)+\dots+(\lambda_p+1)\mu(A_p\cap B) }\Big] \\
&=e^{\psi(\lambda_1)\theta(A_1\setminus B)+\dots+\psi(\lambda_p)\theta(A_p\setminus B)}e^{\psi(\lambda_1+1)\theta(A_1\cap B)+\dots+\psi(\lambda_p+1)\theta(A_p\cap B) } \\
&=e^{\psi(\lambda_1)\theta(A_1 )+\dots+\psi(\lambda_p)\theta(A_p)}e^{(\psi(\lambda_1+1)-\psi(\lambda_1))\theta(A_1\cap B)+\dots+(\psi(\lambda_p+1)-\psi(\lambda_p))\theta(A_p\cap B) } .
\end{align*}
 Then it suffices to notice that:
 $$\psi(q+1)-\psi(q)=m+\sigma^2q+\frac{1}{2} \sigma^2 +\int_\R(e^{(q+1)z}-e^{qz})\nu(dz)$$ and $\psi(1)=0$.\qed

\begin{lemma}\label{eps_gauss}
If the process $\omega_\epsilon$ is defined as $\omega_\epsilon(x) = \mu(A_\epsilon(x))$ where $\mu$ is an independently scattered random measure associated to $(\varphi,\theta)$ with  
$\varphi(q) = -i q\gamma^2/2 - q^2 \gamma^2/2$ and $\theta$ given by \ref{theta}, then:
\begin{equation*}
\lim_{\epsilon \to 0} \E\left[  \frac{e^{\omega_{\epsilon}(x)}}{ z - \int_0^1   K_z(r) e^{\omega_\epsilon(r)} dr } \right] = 
\E\left[ \left(  z - \int_0^1 \left(\frac{\tau}{|r-x|}\right)_+^{\gamma^2} K_z(r) M(dr)  \right)^{-1}\right]
\end{equation*}
where $M$ is the lognormal MRM. 
\end{lemma}
 
\proof

One can check that $(\omega_\epsilon(x))_{x\in [0;1]}$ is a stationary gaussian process with covariance given by $\gamma^2 \rho_\epsilon(x-y)$.
So, using Girsanov transform, we can write:
\begin{align*}
\E\left[  \frac{e^{\omega_{\epsilon}(x)}}{ z - \int_0^1   K_z(r) e^{\omega_\epsilon(r)} dr } \right] = \E\left[  \left( z -\int_0^1  K_z(r) e^{\gamma^2 \rho_\epsilon(r-x)} e^{\omega_\epsilon(r)}  dr \right)^{-1} \right]
\end{align*}

We are interested in the limit when $\epsilon$ goes to $0$ of this latter term, we thus approximate it with a simpler term:
\begin{align}
&\Bigg| \E\left[  \left( z -\int_0^1 K_z(r) e^{\gamma^2 \rho_\epsilon(r-x)} e^{\omega_\epsilon(r)} \,dr \right)^{-1} \right] \nonumber\\
&- \E\left[  \left( z -\int_0^1 K_z(r) \left(\frac{\tau}{|r-x|}\right)_+^{\gamma^2}  e^{\omega_\epsilon(r)}  \,dr \right)^{-1} \right] \Bigg|\nonumber\\
&\leq \frac{1}{|\Im(z)|^2} \E\left[\int_0^1 |K_z(r)| e^{\omega_\epsilon(r)} \left| e^{\gamma^2\rho_\epsilon(r-x)} - \left(\frac{\tau}{|r-x|}\right)_+^{\gamma^2} \right| \,dr \right]
\nonumber\\
&\leq \frac{1}{|\Im(z)|^3} \int_0^1 \left| e^{\gamma^2\rho_\epsilon(r-x)} - \left(\frac{\tau}{|r-x|}\right)_+^{\gamma^2} \right| \,dr \label{exp_rho_eps}
\end{align}
where we have used Lemmas \ref{lem:im} and \ref{imK} and the normalization ${\psi}(1)=0$.

Because $\gamma^2<1$, the dominated convergence theorem implies that \ref{exp_rho_eps} converges to $0$ when $\epsilon$ goes to $0$. 

We thus look at the limit when $\epsilon$ goes to $0$ of the term:
\begin{equation*}
\E\left[  \left( z - \int_0^1 K_z(r) \left(\frac{\tau}{|r-x|}\right)_+^{\gamma^2}  e^{\omega_\epsilon(r)}  \,dr \right)^{-1} \right]. 
\end{equation*}

The random variable 
\begin{equation*}
\int_0^1 K_z(r) \left(\frac{\tau}{|r-x|}\right)_+^{\gamma^2}  M(dr)  
\end{equation*}
is well defined and is finite almost surely since: 
\begin{equation*}
\E\left[\left|\int_0^1 K_z(r) \left(\frac{\tau}{|r-x|}\right)_+^{\gamma^2}  M(dr)\right|\right] \leq \int_0^1 |K_z(r)| \left(\frac{\tau}{|r-x|}\right)_+^{\gamma^2} \,dr < +\infty. 
\end{equation*}

And thus, we can compute:
\begin{align*}
&\Bigg| \E\left[  \left( z -\int_0^1 K_z(r) \left(\frac{\tau}{|r-x|}\right)_+^{\gamma^2}  e^{\omega_\epsilon(r)}  \,dr \right)^{-1} \right] \\
&- \E\left[  \left( z -\int_0^1 K_z(r) \left(\frac{\tau}{|r-x|}\right)_+^{\gamma^2}  M(dr) \right)^{-1} \right]\Bigg|\\
&\leq \frac{1}{|\Im(z)|^2}  \E\left[\left|\int_0^1 K_z(r) \left(\frac{\tau}{|r-x|}\right)_+^{\gamma^2} (e^{\omega_\epsilon(r)} dr -  M(dr))  \right|\right], 
\end{align*}
and, for all $n \in \N$, this latter term is smaller than
\begin{align}
&\E\left[  \left|\int_0^1 K_z(r) \left[ \left(\frac{\tau}{|r-x|}\right)_+^{\gamma^2} - \min\left(\left(\frac{\tau}{|r-x|}\right)_+^{\gamma^2}, n \right) \right] e^{\omega_\epsilon(r)} dr \right|\right] \label{t1} \\
&+ \E\left[  \left|\int_0^1 K_z(r) \min\left(\left(\frac{\tau}{|r-x|}\right)_+^{\gamma^2}, n \right)  (e^{\omega_\epsilon(r)} dr -  M(dr))  \right|\right] \label{t2} \\
&+ \E\left[  \left|\int_0^1 K_z(r) \left[ \left(\frac{\tau}{|r-x|}\right)_+^{\gamma^2} - \min\left(\left(\frac{\tau}{|r-x|}\right)_+^{\gamma^2}, n \right) \right] M(dr) \right|\right].\label{t3}
\end{align}

The two quantities \ref{t1} and \ref{t3} are smaller than
\begin{equation}
\int_0^1 |K_z(r)| \left[ \left(\frac{\tau}{|r-x|}\right)_+^{\gamma^2} - \min\left(\left(\frac{\tau}{|r-x|}\right)_+^{\gamma^2}, n \right) \right] dr 
\end{equation}
and thus converge to $0$, uniformly in $\epsilon$ as $n$ goes to infinity.  

For a fixed $n$, the function $\min((\tau/|r-x|)_+^{\gamma^2}, n)$ is measurable and bounded and thus it is plain to see that, for a fixed $n$, the term \ref{t2} goes to $0$ when $\epsilon$ goes to $0$.

The lemma follows gathering the above estimates. \qed

\begin{lemma}
If the process $\omega_\epsilon$ is defined as $\omega_\epsilon(x) = \mu(A_\epsilon(x))$ where $\mu$ is an independently scattered random measure associated to $(\varphi,\theta)$ where $\varphi$ is given by 
(\ref{char}),i.e. 
\begin{align*}
\varphi(q) = imq - \frac{\gamma^2}{2} q^2 + \int_\R (e^{iqx}-1) \nu(dx)
\end{align*}
and where $\theta$ given by (\ref{theta}), then:
\begin{equation*}
\lim_{\epsilon \to 0} \E\left[  \frac{e^{\omega_{\epsilon}(x)}}{ z - \int_0^1   K_z(r) e^{\omega_\epsilon(r)} dr } \right] = 
\E\left[ \left(  z - \int_0^1 \left(\frac{\tau}{|r-x|}\right)_+^{\kappa} K_z(r) Q(dr)  \right)^{-1}\right]
\end{equation*}
with $\kappa=\gamma^2+\int_\R (e^x-1)^2 \nu(dx)$ and
where the random Radon measure $Q$ is defined, conditionally on a MRM denoted by $M$ whose structure exponent is $\zeta(q):=q-\varphi(-iq)$, as the almost sure weak 
limit as $\epsilon$ goes to $0$ of the family of random measures $Q_\epsilon(dt):=e^{\overline{\omega}_\epsilon(t)} M(dt)$ where, for each $\epsilon>0$, the random process $\overline{\omega}_\epsilon$ is independent of
$M$ and defined as $\overline{\omega}_\epsilon(t)=\overline{\mu}(A_\epsilon(t))$ where $\overline{\mu}$ is the independently scattered log infinitely divisible random measure 
associated to $(\bar{\varphi},\theta(\cdot \cap A_0(x)))$ where 
\begin{equation}
\bar{\varphi}(p) = ip (\gamma^2-\kappa) +\int_\R(e^{ipx}-1)(e^x-1)\nu(dx).
\end{equation}
\end{lemma}
 
\proof
We want to apply Lemma \ref{lem:girs} to the process $\omega_\epsilon$. 
If we set $B=A_\epsilon(x)$, Lemma \ref{lem:girs} tells us that, under $\P_B$, the process $\omega_\epsilon$ possesses the same law as the process 
\begin{equation*}
\omega^{(1)}_\epsilon(r)+\omega^{(2)}_\epsilon(r) \quad \text{with }\omega^{(1)}_\epsilon(r)=\mu^{(1)}(A_\epsilon(r))\text{ and }\omega^{(2)}_\epsilon(r)=\mu^{(2)}(A_\epsilon(r)),
\end{equation*}
where $\mu^{(1)}_\epsilon,\mu^{(2)}_\epsilon$ are independent independently scattered log infinitely divisible random measures respectively associated to $(\varphi,\theta)$ and $(\varphi^{(2)},\theta^{(2)})$ with:
\begin{equation}
\varphi^{(2)}(q) = i\gamma^2q+\int_\R(e^{iqx}-1)(e^x-1)\nu(dx) \text{ and }\theta^{(2)}(\cdot)=\theta(\cdot\cap A_\epsilon(x)).
\end{equation}
Define: 
\begin{align}
\kappa &= \gamma^2 +\int_\R(e^x-1)^2\nu(dx),
\overline{\varphi}(q)=\varphi^{(2)}(q)-iq\kappa,
\overline{\psi}(q)=\overline{\varphi}(-iq).
\end{align}
Notice that $\overline{\psi}$ is then normalized so as to make $\overline{\psi}(1)=\overline{\psi}(0)=0$.
Let us define the process $\overline{\omega}_\epsilon$ by:
\begin{equation}
\overline{\omega}_\epsilon(r)=\omega^{(2)}_\epsilon(r)-\kappa\theta(A_\epsilon(r)\cap A_\epsilon(x))=\omega^{(2)}_\epsilon(r)-\kappa\rho_\epsilon(r-x),
\end{equation}
and notice that $\E[e^{iq \overline{\omega}_\epsilon(r)}]=e^{\overline{\varphi}(q)\rho_\epsilon(r-x)}$.

We can now apply Lemma \ref{lem:girs}:
\begin{align*}
\E\left[  \frac{e^{\omega_{\epsilon}(x)}}{ z -\int_0^1   K_z(r) e^{\omega_\epsilon(r)} dr } \right] = \E\left[  \left( z -\int_0^1 K_z(r) e^{\overline{\omega}_\epsilon(r)+\kappa\rho_\epsilon(r-x)+\omega_\epsilon(r)} 
 dr \right)^{-1} \right]  
\end{align*}


We are interested in the limit when $\epsilon$ goes to $0$ of this latter term, we thus approximate it with a simpler term:
\begin{align}
&\Bigg| \E\left[  \left( z -\int_0^1  e^{\overline{\omega}_\epsilon(r)+\omega_\epsilon(r)+\kappa\rho_\epsilon(r-x)} K_z(r) \,dr \right)^{-1} \right] \nonumber\\
&- \E\left[  \left( z -\int_0^1  e^{\overline{\omega}_\epsilon(r)+\omega_\epsilon(r)} \left(\frac{\tau}{|r-x|}\right)_+^{\kappa} K_z(r) \,dr \right)^{-1} \right] \Bigg|\nonumber\\
&\leq \frac{1}{|\Im(z)|^2} \E\left[\int_0^1 e^{\overline{\omega}_\epsilon(r)+\omega_\epsilon(r)} \left| e^{\kappa\rho_\epsilon(r-x)} - \left(\frac{\tau}{|r-x|}\right)_+^{\kappa} \right| |K_z(r)| \,dr \right]
\nonumber\\
&\leq \frac{1}{|\Im(z)|^3} \int_0^1 \left| e^{\kappa\rho_\epsilon(r-x)} - \left(\frac{\tau}{|r-x|}\right)_+^{\kappa} \right| \,dr \label{exp_rho_eps_2}
\end{align}
where we have used Lemmas \ref{lem:im} and \ref{imK}, the normalizations $\overline{\psi}(1)=0, {\psi}(1)=0$ and the independence between $\overline{\omega}_\epsilon$ and $\omega_\epsilon$.

Let us show that $\kappa <1$. Indeed, we have:
\begin{align*} 
\kappa &= \gamma^2 +\int_\R(e^x-1)^2\nu(dx) \\
&= \gamma^2 + \int_\R (e^{2x}-1) \nu(dx) - 2 \int_\R (e^{x}-1) \nu(dx) \\
&= \gamma^2 + \int_\R (e^{2x}-1) \nu(dx) + 2 (m+\frac{1}{2}\gamma^2) \\
&= 2m + 2 \gamma^2 +\int_\R (e^{2x}-1) \nu(dx) \\
&= \psi(2)
\end{align*}
where, in the third line, we used the fact that $\psi(1) = 0$ (which implies the relation $\int_\R (e^{x}-1) \nu(dx) = -(m+\gamma^2/2)$). We will now show that $\psi(2)$ is strictly less than $1$.
It suffices to show that $\zeta(2)>1$. Using the concavity of the function $\zeta$, we have the inequality:
\begin{equation}
\frac{\zeta(2+\epsilon) - \zeta(1)}{1+\epsilon} < \zeta(2) - \zeta(1)
\end{equation}
and with assumption \ref{condzeta}, we see that $\zeta(2) - \zeta(1) = \zeta(2) - 1 >0$. We can thus conclude that $\kappa <1$. 
 
Because $\kappa<1$, the dominated convergence theorem implies that \ref{exp_rho_eps_2} converges to $0$ when $\epsilon$ goes to $0$.

For each Borelian set $A$ of $[0;1]$, the family $M_\epsilon(A) := \int_A e^{\omega_\epsilon(r)} dr, \epsilon >0$ is a positive martingale with respect to $\epsilon$ and that it converges almost surely to $M(A)$. 
With the assumption $\ref{condzeta}$ and in particular the condition $\zeta(2+\epsilon)>1$, we can show (see \cite{Bacry} for a proof) that the family $(M_\epsilon(A))_{\epsilon >0}$ is in fact uniformly integrable. 
In particular, if we let $\mathcal{F}_{\epsilon}$ be the sigma field generated by the family of random variables $(\omega_\eta(r))_{\eta >\epsilon, r\in \R}$, we have the following almost sure equality:
\begin{equation}
\E\left[M(A) |\mathcal{F}_\epsilon\right] = M_\epsilon(A).
\end{equation}

Conditionally to the random measure $M$, the family $P_\epsilon(A):= \int_A e^{\overline{\omega}_\epsilon(r)} M(dr), \epsilon >0$ is also a positive martingale with respect to $\epsilon$. Thus, $P_\epsilon(A)$ converges 
almost surely to a random variable
that we will denote by $P(A)$. We know that this defines a random Radon measure $P$ on $[0;1]$ and that the family of random Radon measures $P_\epsilon$ converges, when $\epsilon$ goes to $0$, weakly almost surely to 
$P$ in the space of Radon measures. 
Denote, conditionally to the random measure $M$, by $\P_M$ the law $\P[\cdot|M]$ and let us show that the family $(P_\epsilon([0;1]))_{\epsilon >0}$ is $\P_M$-uniformly integrable. 
Let $\overline{\delta}$ be such that $\overline{\psi}(1+\overline{\delta}) < +\infty$
(we can show, using the condition $\psi(2+\delta)<+\infty$, that that there exists such $\overline{\delta}$ ). 
We will show that the family $(P_\epsilon([0;1]))_{\epsilon >0}$ is uniformly bounded in $L^{1+\overline{\delta}}(\P_M)$.
Indeed, conditionally to the random measure $M$:
\begin{align*}
\E_M\left[\left(\int_0^1 e^{\overline{\omega}_\epsilon(r)} M(dr)\right)^{1+\overline{\delta}} \right] &
\leq \E_M\left[ \int_0^1 e^{(1+\overline{\delta})\overline{\omega}_\epsilon(r)} M(dr) \right] M[0;1]^{\overline{\delta}} \\
&\leq \int_0^1 e^{\overline{\psi}(1+\overline{\delta})\rho_\epsilon(r-x)} M(dr) M[0;1]^{\overline{\delta}} \\ 
& \leq M[0;1]^{\overline{\delta}}  e^{\overline{\psi}(1+\overline{\delta})} \int_0^1  \left(\frac{\tau}{|r-x|}\right)_+^{\kappa} M(dr) < +\infty.
\end{align*} 
The family $(P_\epsilon([0;1]))_{\epsilon >0}$ is therefore $\P_M$-uniformly integrable, in particular, $P_\epsilon([0;1])$ converges to $P([0;1])$ also in $L^1$, which implies that $P$ is a non degenerated 
random measure. 
Moreover, denoting 
by $\overline{\mathcal{F}}_{\epsilon}$ the sigma field generated by the family of random variables $(\overline{\omega}_\eta(r))_{\eta >\epsilon, r\in \R}$, 
we have, almost surely, conditionally to $M$, for all Borelian set $A$ of $[0;1]$:
\begin{equation*}
\E_M\left[P(A) |\mathcal{F}_\epsilon\right] = P_\epsilon(A).
\end{equation*}

Now, as before, it is easy to see that the family $Q_\epsilon(A):=\int_A e^{\omega_\epsilon(r)+\overline{\omega}_\epsilon(r)} dr, \epsilon >0$ is also a positive martingale with respect to $\epsilon$. Therefore, 
$Q_\epsilon(A)$ converges almost surely to a random variable that we will denote by $Q(A)$. This defines a random Radon measure $Q$ and the family of random Radon measure $Q_\epsilon$
converges, as $\epsilon \to 0$, weakly almost surely to $Q$ in the space of Radon measure. We want to show that the two random measures $P$ and $Q$ have the same law. 

Gathering the above arguments, we can write, almost surely:
\begin{align*}
\E\left[P(A)|\sigma(\mathcal{F}_\epsilon,\overline{\mathcal{F}}_\epsilon)\right] &= \E\left[\E[P(A)|\overline{\mathcal{F}}_\epsilon] \right] \\
&= \E\left[\int_A e^{\overline{\omega}_\epsilon(r)} M(dr) |\mathcal{F}_\epsilon\right] \\
&= \int_A e^{\omega_\epsilon(r)+\overline{\omega}_\epsilon(r)} dr, 
\end{align*} 
and the latter quantity has the same law as $Q_\epsilon(A)$. Since the martingale $(\E[P(A)|\sigma(\mathcal{F}_\epsilon,\overline{\mathcal{F}}_\epsilon)])_{\epsilon>0}$ is uniformly integrable, we deduce that the 
family $(Q_\epsilon(A))_{\epsilon>0}$ is also uniformly integrable. Hence, both random variables $P(A)$ and $Q(A)$ have the same law. We can show easily that in fact the two random measures $P$ and $Q$
have the same law. In particular, $Q$ is non degenerated. 

It is now easy to see that, for all bounded and continuous function $f$, the two random variables $\int_\R f(r) P(dr)$ and $\int_\R f(r) Q(dr)$ have the same law.  
By regularizing the function $\left(\frac{\tau}{|r-x|}\right)_+^{\kappa}$ and with the dominated convergence theorem, we conclude as in the proof of lemma \ref{eps_gauss} using the fact that $\kappa<1$ that:
\begin{equation}
\int_0^1 K_z(r) \left(\frac{\tau}{|r-x|}\right)_+^{\kappa}  Q(dr) \stackrel{(law)}{=} \int_0^1 K_z(r) \left(\frac{\tau}{|r-x|}\right)_+^{\kappa}  P(dr). 
\end{equation}
Gathering the above argument and letting $\epsilon$ go to $0$ concludes the proof. \qed


\end{document}